\documentclass[10pt]{amsart}

\usepackage{hyperref}
\usepackage{enumerate}
\input{diagrams}

\newtheorem{proposition}{Proposition}[section]
\newtheorem{lemma}[proposition]{Lemma}
\newtheorem{corollary}[proposition]{Corollary}
\newtheorem{theorem}[proposition]{Theorem}

\theoremstyle{definition}

\newtheorem{example}[proposition]{Example}

\theoremstyle{remark}
\newtheorem{remark}[proposition]{Remark}

\newcommand{\thlabel}[1]{\label{th:#1}}
\newcommand{\thref}[1]{Theorem~\ref{th:#1}}
\newcommand{\selabel}[1]{\label{se:#1}}
\newcommand{\seref}[1]{Section~\ref{se:#1}}
\newcommand{\lelabel}[1]{\label{le:#1}}
\newcommand{\leref}[1]{Lemma~\ref{le:#1}}
\newcommand{\prlabel}[1]{\label{pr:#1}}
\newcommand{\prref}[1]{Proposition~\ref{pr:#1}}
\newcommand{\colabel}[1]{\label{co:#1}}

\newcommand{\exlabel}[1]{\label{ex:#1}}
\newcommand{\exref}[1]{Example~\ref{ex:#1}}

\newcommand{\eqlabel}[1]{\label{eq:#1}}
\newcommand{\equref}[1]{(\ref{eq:#1})}

\newcommand{\Hom}{{\rm Hom}}
\newcommand{\End}{{\rm End}}

\newcommand{\Ker}{{\rm Ker}\,}

\newcommand{\can}{{\rm can}}

\def\Ab{\underline{\underline{\rm Ab}}}

\def\ot{\otimes}
\def\ev{{\rm ev}}
\def\coev{{\rm coev}}

\def\colim{{\rm colim}\,}

\newcommand{\Aa}{\mathcal{A}}
\newcommand{\Bb}{\mathcal{B}}
\newcommand{\Cc}{\mathcal{C}}
\newcommand{\Dd}{\mathcal{D}}

\newcommand{\Ff}{\mathcal{F}}
\newcommand{\Gg}{\mathcal{G}}

\newcommand{\Mm}{\mathcal{M}}

\newcommand{\Yy}{{\mathcal Y}}
\newcommand{\Zz}{{\mathcal Z}}

\def\*C{{}^*\hspace*{-1pt}{\Aa}}
\def\text#1{{\rm {\rm #1}}}

\def\ul{\underline}

\def\equal#1{\smash{\mathop{=}\limits^{#1}}}

\begin{document}
\title[Constructing infinite comatrix corings
from colimits ]{Constructing infinite comatrix corings
from colimits}
\author{S. Caenepeel}
\address{Faculty of Engineering,
Vrije Universiteit Brussel, VUB, B-1050 Brussels, Belgium}
\email{scaenepe@vub.ac.be}
\urladdr{http://homepages.vub.ac.be/\~{}scaenepe/}
\author{E. De Groot}
\address{Faculty of Engineering,
Vrije Universiteit Brussel, VUB, B-1050 Brussels, Belgium}
\email{edegroot@vub.ac.be}
\urladdr{http://homepages.vub.ac.be/\~{}edegroot/}
\author{J. Vercruysse}
\address{Faculty of Engineering,
Vrije Universiteit Brussel, VUB, B-1050 Brussels, Belgium}
\email{joost.vercruysse@vub.ac.be}
\urladdr{http://homepages.vub.ac.be/\~{}jvercruy/}
\thanks{}
\subjclass{16W30}

\keywords{}

\begin{abstract}
We propose a class of infinite comatrix corings, and describe them as
colimits of systems of usual comatrix corings. The infinite comatrix corings
of El Kaoutit and G\'omez Torrecillas are special cases of our construction,
which in turn can be considered as a special case of the comatrix corings
introduced recently by G\'omez Torrecillas an the third author.
\end{abstract}

\maketitle

\section*{Introduction}
Corings were introduced by Sweedler in 1975 \cite{Sweedler65}; since the
beginning of the century, there has been a renewed interest in corings,
initiated by an observation made by Takeuchi that most type of modules that
are considered in Hopf algebra theory, like Hopf modules, Yetter-Drinfeld
modules, entwined modules, are in fact comodules over certain corings.
A detailed discussion of recent applications of corings can be found in
\cite{BrzezinskiWisbauer}.\\
One of the beautiful applications is a reformulation of descent theory and
Galois theory. To a ring morphism $B\to A$, we can associate a
coring $A\ot_B A$, called Sweedler's canonical coring, and the category of descent data
is isomorphic to the category of comodules over the coring. To an action or
coaction of a group or Hopf algebra on $A$, we can associate a coring, and
there exists a canonical coring map from Sweedler's coring to this coring.
A necessary condition for the Galois descent is that this map is an isomorphism.
This was observed by Brzezi\'nski in his paper \cite{Brzezinski02}, see also
\cite{Caenepeel03} for a detailed discussion.\\
A more general theory was proposed by El Kaoutit and G\'omez Torrecillas
\cite{Kaoutit}. We start from two rings $A$ and $B$, connected by a $(B,A)$-module
$P$. If $P$ is finitely generated and projective as a right $A$-module,
$P^*\ot_B P$ is an $A$-coring. If $A$ and $B$ are connected via a ring morphism
$B\to A$, then we can take $P=A$ considered as a $(B,A)$-bimodule, and we
recover Sweedler's coring. $P^*\ot_B P$ is called a comatrix coring, and
several properties of the theory outlined in \cite{Caenepeel03} can be generalized,
we refer to \cite{CaenepeelDGV} and \cite{Kaoutit}.\\
The condition that $P$ is finitely generated and projective as a right $A$-module
is crucial in the theory. Nevertheless, El Kaoutit and G\'omez Torrecillas \cite{KGT}
proposed an infinite version of comatrix corings, starting from an infinite
collection of finitely generated projective right $A$-modules $\{P_i~|~i\in I\}$.
They consider the direct sum $P$ of the $P_i$, and the direct sum $P^\dagger$
of the $P_i^*$. The tensor product of $P^\dagger$ and $P$ over a suitable ring
$R$ is then a coring, called the infinite comatrix coring. They give several
descriptions and properties of this coring, including a version of the Faithfully
Flat Descent Theorem. One of the important features is the fact that the
ring $R$ has no unit; it is a ring with orthogonal idempotent local units.\\
The natural framework needed to introduce infinite comatrix corings was
proposed recently by G\'omez Torrecillas and the third author in \cite{GTV}.
The philosophy is the following. Let $P$ be a $(B,A)$-bimodule, and consider the
functor $F=-\ot_B P:\ \Mm_B\to \Mm_A$. $F$ has a right adjoint $G$ of the form
$-\ot_A Q$ for some $(B,A)$-bimodule $Q$ if and only if $P$ is finitely generated
and projective a right $A$-module, and in this case $Q\cong P^*$. Then we have
a so-called comatrix coring context (see \cite{BrzezinskiG} or \cite{Caenepeel98}
for the definition). In fact we need such an adjunction to be able to define the
coproduct on the comatrix coring. Instead of considering rings with a unit, we
now consider firm rings, these are rings $A$ having the property that
the canonical map $A\ot_AA\to A$ is an isomorphism. Firm bimodules over firm rings
form a monoidal category, and we can consider corings over firm rings. 
The bimodules in a comatrix coring context connecting firm rings
are not necessarily finitely generated projective. The comatrix coring contexts
from \cite{KGT} are of this type: only one of the two rings involved has units,
the other one has only local units (a complete set of orthogonal idempotents).\\
In this paper, we propose some classes of infinite comatrix corings. In the first three
Sections, we have collected some necessary preliminary results: in
\seref{1}, we briefly introduce the comatrix corings from \cite{GTV},
and recall some of the elementary results, for example the Faithfully
Flat Descent Theorem; in \seref{2}
we show that the colimit of a functor that is a coalgebra in a functor category
is itself a coalgebra; in \seref{3}, we discuss split directed systems and their
colimits. The main results appear in \seref{4}: we describe comatrix corings
associated associated to an $(A,B)$-bimodule
$P$, where $A$ is a ring with unit, and $B$ a ring with idempotent local units. These rings are
 colimits of a (split) directed system of rings with unit and 
we can describe the category of firm $B$-modules.
The comatrix
corings can also be described as colimits. 
Firm modules over
rings with idempotent local units are constructed in \seref{6}: we consider a split direct system $\ul{M}$
in some $k$-linear category $\Aa$ with a colimit, and a product preserving
functor $\omega$ to the category of right $A$-modules. The ring $B$ is then the
colimit of the $\Aa$-endomorphism rings of the $M_i$, and the $\omega(M_i)=P_i$
form a split direct system of $(B_i,A)$-bimodules. An interesting special case is
considered in \seref{7}: we consider an $A$-coring $\Cc$, and let $\Aa=\Mm^\Cc$,
and $\omega$ the functor forgetting the coaction. In this situation, we can define
a canonical coring map from the associated comatrix coring to $\Dd$, and $\ul{M}$
is called a system of Galois $\Cc$-comodules if this map is an isomorphism. 
The comatrix corings introduced in \cite{KGT} are special cases.

\section{Comatrix corings over firm algebras}\selabel{1}
Let $k$ be a commutative ring, and $A$ a $k$-algebra, not necessarily with a unit.
If $A$ has a unit, then the canonical map $A\ot_A A\to A$ is an isomorphism,
but not conversely. We say that $A$ is a {\sl firm algebra} if $A\ot_A A\to A$
is an isomorphism. In \cite{Taylor}, firm algebras are called {\sl regular algebras};
in \cite{Caenepeel98}, they are called {\sl unital}.
Algebras with local units are firm.\\
Let $A$ be a firm algebra. A right $A$-module is called {\sl firm}
if the canonical map $M\ot_A A\to M$ is an isomorphism. If $A$ is an algebra with
unit, then all modules are firm. $\Mm_A$ will be the category of {\sl firm} right $A$-modules
and right $A$-linear maps. In a similar way, we introduce the categories of firm
left modules and firm bimodules. The category of all (not nessecary firm) right $A$-modules will be denoted by $\widetilde{\Mm}_A$.  A (firm) left $A$-module is called flat if the functor $-\ot_AM:\widetilde{\Mm}_A\to \Ab$ is exact. The category $\Mm_A$ of firm right $A$-modules over a firm ring $A$ is always an abelian category. Under the extra condition that $A$ is flat as left $A$-module, the kernels in $\Mm_A$ can be computed already in $\Ab$. If $A$ is a ring with local units, i.e. for every $a\in A$, there exists an $e\in A$ such that $ae=ea=a$, then $A$ is a firm ring and $A$ is flat as a left and right $A$-module ($A$ is even locally projective as a left and right $A$-module), so in this situation the kernels of $\Mm_A$ and ${_A\Mm}$ can be computed in $\Ab$.\\
The category ${}_A\Mm_A$ of firm $A$-bimodules is
a monoidal category, so we can consider corings over firm algebras, these are
coalgebras in the monoidal category ${}_A\Mm_A$. If $\Cc$ is a coring over a firm
$k$-algebra $A$, then we can define left and right $\Cc$-comodules. A right $\Cc$-comodule
$(M,\rho^r)$ is a firm right $A$-module together with a right $A$-linear map
$\rho^r:\ M\to M\ot_A\Cc$ satisfying the usual coassociativity and counit properties. The category of right $\Cc$-comodules and $\Cc$-colinear maps is denoted by $\Mm^\Cc$.
Similary one introduces categories ${^\Cc\Mm}$, ${_B\Mm^\Cc}$ and ${^\Dd\Mm^\Cc}$. For $M\in{}^\Cc\Mm$ and $N\in\Mm^\Cc$, we define the cotensor product $N\ot^\Cc M$
as the following equalizer in $\Ab$:
\begin{equation}\eqlabel{1}
N\ot^\Cc M \rTo^{} N\ot_AM
\pile{\rTo^{N\ot_A\rho^l}\\ \rTo_{\rho^r\ot_AM}}
N\ot_A\Cc\ot_AM.
\end{equation}
If $M\in {}^\Cc\Mm_B$, where $B$ is a firm ring, which is flat as a left $B$-module, then
\equref{1} is also an equalizer in $\Mm_B$, hence $N\ot^\Cc M$ is a firm right $B$-module.\\
Let $A$ and $B$ be firm $k$-algebras. The notion of Morita context can be generalized,
by requiring that the connecting bimodules are firm. If one of the morphisms in a Morita
context is bijective, then we can associate a pair of adjoint functors to it. More generally,
we have the following result (see \cite[Theorem 1.1.3]{Caenepeel98}).

\begin{proposition}\prlabel{1.1}
Let $B$ and $A$ be firm $k$-algebras, and $P\in {}_B\Mm_A$,
$P^{\dagger}\in {}_A\Mm_B$ firm bimodules. Consider two bimodule maps
$$\eta:\ B\to P\ot_AP^{\dagger}~~{\rm and}~~
\varepsilon:\ P^{\dagger}\ot_BP\to A.$$
We use the following Sweedler-type notation:
$$\eta(b)=b^-\ot_Ab^+\in P\ot_AP^{\dagger},$$
where summation is implicitly understood, as usual. Assume that $\eta$ and $\varepsilon$
satisfy the following formulas, for all $b\in B$, $p\in P$, $q\in P^{\dagger}$:
\begin{equation}\eqlabel{1.1.1}
b^-\varepsilon(b^+\ot_Ap)=bp~~;~~\varepsilon(q\ot_A b^-)b^+=qb.
\end{equation}
Then we have a pair of adjoint functors $(F,G)$
$$F=-\ot_BP:\ \Mm_B\to \Mm_A~~;~~G=-\ot_AP^\dagger:\
\Mm_A\to \Mm_B.$$
\end{proposition}

\begin{proof}
The unit and counit of the adjunction are 
$$\eta_N=N\ot_B \eta:\ N\to N\ot_BP\ot_AP^\dagger~~;~~
\varepsilon_M=M\ot_A\varepsilon:\ M\ot_AP^\dagger\ot_BP\to M,$$
for all $N\in \Mm_B$, $M\in \Mm_A$.
\end{proof}

Following \cite{BrzezinskiG},
$(B,A,P,P^\dagger,\eta,\varepsilon)$ is called a 
{\sl comatrix coring context}\index{comatrix coring context}. 
To a comatrix coring context $(B,A,P,P^\dagger,\eta,\varepsilon)$, we can 
associate an $A$-coring $\Dd$ (called comatrix coring) and a $B$-ring $\Aa$
(called matrix ring, or elementary algebra (see \cite{Caenepeel98,Taylor}).
They are given by the following data:\\
$\Dd=P^\dagger\ot_BP$, with
$$\Delta_\Dd=P^\dagger\ot_B\eta\ot_BP:\ \Dd\to \Dd\ot_A\Dd~~;~~
\varepsilon_\Dd=\varepsilon:\ \Dd\to A.$$
$\Aa=P\ot_AP^\dagger$, with
$$m_\Aa=P\ot_A\varepsilon\ot_AP^\dagger:\ \Aa\ot_B\Aa\to \Aa~~;~~
\eta_\Aa=\eta:\ B\to P\ot_AP^\dagger.$$
$P$ is a right $\Dd$-comodule, and $P^\dagger$ is a left $\Dd$-comodule;
the right and left coactions are the following:
$$\rho^r:\ P\to P\ot_AP^\dagger\ot_BP,~~
\rho^r(bp)=b^-\ot_Ab^+\ot_B p$$
$$\rho^l:\ P^\dagger\to P^\dagger\ot_BP\ot_AP^\dagger,~~
\rho^l(qb)=q\ot_B b^-\ot_Ab^+$$

\begin{proposition}\prlabel{1.2}
Let $(B,A,P,P^\dagger,\eta,\varepsilon)$ be a comatrix coring context, and
assume that $B$ is flat as left $B$-module.
Then we have
a pair of adjoint functors $(K,R)$
$$K=-\ot_BP:\ \Mm_B\to \Mm^\Dd~~;~~
R=-\ot^\Dd P^\dagger:\ \Mm^\Dd\to \Mm_B.$$
\end{proposition}

\begin{proof}
It follows from the comments preceeding \prref{1.1} that $R(M)$ is a firm right
$B$-module, for every $M\in \Mm^\Dd$.
We restrict to giving the unit and the counit of the adjunction. For $N\in \Mm_B$
and $M\in \Mm^\Dd$, we have
$$\eta_N:\ N\to (N\ot_BP)\ot^\Dd P^\dagger,~~
\eta_N(nb)=n\ot_B b^-\ot_A b^+;$$
$$\varepsilon_M:\ (M\ot^\Dd P^\dagger)\ot_BP\to M,~~
\varepsilon_M(\sum_j m_j\ot_Aq_j\ot p_j)=\sum_jm_j\varepsilon(q_j\ot_Bp_j).$$
Let us show that $\eta_N(n)\in (N\ot_BP)\ot^\Dd P^\dagger$,
for all $n\in N$. Since $N$ is firm as a right $B$-module, it suffices to look
at elements of the form $nbcd$, with $n\in N$, $b,c,d\in B$. Since $\eta$ is a
$B$-bimodule map, we have, for all $b,c\in B$ that
$\eta(bc)=b\eta(c)=\eta(b)c$, or
\begin{equation}\eqlabel{1.2.1}
(bc)^-\ot_A (bc)^+=bc^-\ot_A c^+= b^-\ot_Ab^+c.
\end{equation}
Using \equref{1.2.1}, we find easily that
\begin{eqnarray*}
&&\hspace*{-2cm}
(\rho^r_{N\ot_BP}\ot_A P^\dagger)(\eta_N(nbcd))=
(\rho^r_{N\ot_BP}\ot_AP^\dagger)(nb\ot_B c^-\ot_Ac^+d)\\
&=& (N\ot\rho^r\ot P^{\dagger})(n\ot_B bcd^-\ot_Ad^+)\\
&=& n\ot_B b^-\ot_A b^+\ot_B cd^- \ot_A d^+\\
&=& n\ot_B b^-\ot_A b^+c\ot_B d^- \ot_A d^+\\
&=& nb\ot_B c^-\ot_A c^+\ot_B d^- \ot_A d^+\\
&=& ((N\ot_B P)\ot_A\rho^l)(nb\ot_B c^-\ot_A c^+d)\\
&=& ((N\ot_B P)\ot_A\rho^l)(\eta_N(nbcd)).
\end{eqnarray*}
\end{proof}

\begin{theorem}\thlabel{1.3} {\bf (Faithfully flat descent)}
Let $(B,A,P,P^\dagger,f,g)$ be a comatrix coring context, and assume that
$B$ and
$P$ are flat as a left $B$-module. Then $R$ is fully faithful.
$(K,R)$ is a pair of inverse equivalences if and only if $P$ is
faithfully flat as a left $B$-module. 
\end{theorem}

\begin{proof}
Take $M\in \Mm^\Dd$.
If $P\in {}_B \Mm$ is flat, then the map
$$j:\ (M\ot^\Dd P^\dagger)\ot_BP\to M\ot^\Dd (P^\dagger\ot_BP),~~
j((\sum_i m_i\ot_A p_i)\ot_B q)=\sum_i m_i\ot_A (p_i\ot_Bq)$$
is an isomorphism. The map
$$M\ot_A \varepsilon:\ M\ot^\Dd (P^\dagger\ot_BP)\to M\ot_AA\cong M$$
is an isomorphism. If $P\in  {}_B\Mm$ is flat, then $\varepsilon_M=
(M\ot_A \varepsilon)\circ j$ is an isomorphism.\\
Assume that $P\in {}_B\Mm$ is flat. We have to show that
$\eta_N$ is an isomorphism, for every $N\in \Mm_B$. It suffices to show that
the sequence
$$S:~~
0\rTo^{}N\rTo^{N\ot_B \eta}N\ot_BP\ot_AP^\dagger
\pile{\rTo^{N\ot\rho^r\ot P^\dagger}\\ \rTo_{N\ot P\ot \rho^l}}
N\ot_BP\ot_AP^\dagger\ot_BP\ot_AP^\dagger$$
is exact. Since $P\in {}_B\Mm$ is faithfully flat, it suffices to show that
$S\ot_BP$ is exact. It is clear that the sequence is a complex.\\
We first show that $N\ot_B\eta\ot_BP$ is injective: if
$$0=(N\ot_B\eta\ot_BP)(\sum_j n_jb_j\ot_B p_j)=
\sum_j n_j\ot_B b_j^-\ot_A b_j^+\ot_B p_j,$$
then
$$0=\sum_j n_j\ot_B b_j^-\varepsilon( b_j^+\ot_B p_j)=\sum_jn_j\ot_B b_jp_j.$$
Now assume that
$$x=\sum_j n_j\ot_B b_jp_j\ot_A q_j\ot_B c_jr_j\in
\Ker\bigl(N\ot_BP\ot_A\rho^l\ot_BP-
N\ot_B\rho^r\ot_AP^\dagger\ot_BP\bigr).$$
Then
$$\sum_j n_j\ot_B b_jp_j\ot_A q_j\ot_Bc_j^-\ot_Ac_j^+\ot_B r_j=
\sum_j n_j\ot_B b_j^-\ot_Ab_j^+\ot_Bp_j\ot_A q_jc_j\ot_B r_j,$$
and it follows that
\begin{eqnarray*}
x&=&\sum_j n_j\ot_B b_jp_j\ot_A q_j\ot_Bc_j^-\varepsilon(c_j^+\ot_B r_j)\\
&=& \sum_j n_j\ot_B b_j^-\ot_Ab_j^+\ot_Bp_j\varepsilon(q_jc_j\ot_B r_j)\\
&=& (N\ot_B \eta\ot_BP)\bigl(\sum_j n_jb_j\ot_Bp_j\varepsilon(q_jc_j\ot_B r_j).
\end{eqnarray*}
\end{proof}

\section{Corings from colimits}\selabel{2}
Let $F:\ \Zz\to \Mm$ be a covariant functor. Recall (see for example \cite{Borceux}) that a
{\sl cocone}\index{cocone} on $F$ is a couple $(M,m)$ where $M\in \Mm$ and
$m_Z:\ F(Z)\to M$ is a morphism in $\Mm$, for every $Z\in \Zz$, such that
\begin{equation}\eqlabel{2.1.1}
m_Z=m_{Z'}\circ F(f),
\end{equation}
for every $f:\ Z\to Z'$ in $\Zz$. The {\sl colimit}\index{colimit} of $F$ is a cocone
$(C,c)$ on $F$ satisfying the following universal property: if $(M,m)$ is a cocone on $F$,
then there exists a unique morphism $f:\ C\to M$ in $\Mm$ such that
\begin{equation}\eqlabel{2.1.2}
f\circ c_Z=m_Z,
\end{equation}
for every $Z\in \Zz$. If the colimit exists, then it is unique up to isomorphism. We then write
$\colim F=\colim F(Z)=(C,c)$.\\
The colimit $(C,c)$ has the following property: if $f,g:\ C\to M$ are two morphisms in $\Mm$
such that $f\circ c_Z=g\circ c_Z$, for all $Z\in \Zz$, then $f=g$. Indeed, $(M,f\circ c=g\circ c)$
is a cocone on $F$, and $f=g$ follows from the uniqueness in the definition of colimit.\\

From now on, let $\Zz$ be a (small) category and let $(\Mm,\ot, A)$ be a monoidal category. Then $({\rm Func}(\Zz,\Mm),\ot, A)$ is
also a monoidal category. The tensor $\ot$ and the unit $A$ are given by the following formulas:
$$(F\ot G)(Z)=F(Z)\ot G(Z)~~~{\rm and}~~(F\ot G)(f)=F(f)\ot G(f);$$
$$A(Z)=A~~{\rm and}~~A(f)=A,$$
for all $F,G:\ \Zz\to \Mm$,
$Z,Z'\in \Zz$ and $f:\ Z\to Z'$ in $\Zz$. A coalgebra in $({\rm Func}(\Zz,\Mm),\ot, A)$
will be called a $\Zz$-coalgebra in $\Mm$. The result of this Section is the following.

\begin{proposition}\prlabel{2.1}
Let $(G,\Delta,\varepsilon)$ be a $\Zz$-coalgebra in $\Mm$, and assume that
$\colim G=(C,c)$ exists. Then $C$ is a coalgebra in $\Mm$.
\end{proposition}

\begin{proof}
We give a proof of the statement in case of a strict monoidal category $\Mm$. Recall that this is no restriction since every monoidal category is equivalent to a strict monoidal category,
see for example \cite[Prop. IX.5.1]{Kassel}.

For every $Z\in \Zz$, consider the morphism
$$d_Z=(c_Z\ot c_Z)\circ\Delta_Z:\ G(Z)\to C\ot C.$$
Let $f:\ Z\to Z'$ be a morphism in $\Zz$, and look at the diagram
$$\begin{diagram}
G(Z)&\rTo^{\Delta_Z}&G(Z)\ot G(Z)& \rTo^{c_Z\ot c_Z} &C\ot C\\
\dTo^{G(f)}&&\dTo_{G(f)\ot G(f)}&&\dTo_{=}\\
G(Z')&\rTo^{\Delta_{Z'}}&G(Z')\ot G(Z')& \rTo^{c_{Z'}\ot c_{Z'}} &C\ot C
\end{diagram}$$
The left hand square commutes since $\Delta:\ G\to G\ot G$ is a natural transformation,
and the right hand square commutes because $(C,c)$ is a cocone on $G$. It
follows that $(C\ot C,d)$ is a cocone on $G$, and we conclude that there exists a morphism
$\Delta_C:\ C\to C\ot C$
in $\Mm$ such that
$$\Delta_C\circ c_Z=d_Z=(c_Z\ot c_Z)\circ \Delta_Z,$$
for all $Z\in \Zz$. We then have
\begin{eqnarray*}
&&\hspace*{-2cm}
(\Delta_C\ot C)\circ\Delta_C\circ c_Z=(\Delta_C\ot C)\circ (c_Z\ot c_Z)\circ\Delta_Z\\
&=& (c_Z\ot c_Z\ot C)\circ (\Delta_Z\ot c_Z)\circ \Delta_Z\\
&=& (c_Z\ot c_Z\ot c_Z)\circ (\Delta_Z\ot G(Z))\circ \Delta_Z\\
&=& (c_Z\ot c_Z\ot c_Z)\circ (G(Z)\ot \Delta_Z)\circ \Delta_Z\\
&=&
(C\ot \Delta_C)\circ\Delta_C\circ c_Z,
\end{eqnarray*}
for all $Z\in \Zz$. It follows (see \cite[Prop. 2.6.4]{Borceux}) that
$(\Delta_C\ot C)\circ\Delta_C=(C\ot \Delta_C)\circ\Delta_C$, so $\Delta_C$ is a
coassociative comultiplication on $C$.\\
The counit is defined in a similar way: $(A,\varepsilon)$ is a cocone on $G$, so there
exists a morphism $\varepsilon_C:\ C\to A$ in $\Mm$ such that $\varepsilon_C\circ
c_Z=\varepsilon_Z$, for all $Z\in \Zz$. The counit property is verified as follows: for
all $Z\in \Zz$, we have
$$
(\varepsilon_C\ot C)\circ\Delta_C\circ c_Z=
(\varepsilon_C\ot C)\circ(c_Z\ot c_Z)\circ \Delta_Z=
(A\ot c_Z)\circ (\varepsilon_Z\ot G(Z))\circ\Delta_Z=c_Z.
$$
\end{proof}

\begin{proposition}\prlabel{2.2}
Let $(G,\Delta,\varepsilon)$ be as in \prref{2.1}. If $(H,\rho)$ is a right $G$-comodule,
and $\colim H=(M,m)$ exists, then $M$ is a right $C$-comodule.
\end{proposition}

\begin{proof}
For every $Z\in \Zz$, consider the composition
$$r_Z=(m_Z\ot c_Z)\circ \rho_Z:\ H(Z)\to H(Z)\ot G(Z)\to M\ot C.$$
Arguments similar to the ones presented above show that $(M\ot C,r)$ is a cocone on $H$.
It follows that there exists a morphism $\rho_M:\ M\to M\ot C$ such that
$\rho_M\circ m_Z=r_Z$, for every $Z\in \Zz$. Standard computations show that
$\rho_M$ is coassociative and satisfies the counit property.
\end{proof}

\section{Split direct systems}\selabel{3}
Recall that a partially ordered set $(I,\leq)$ is called {\sl directed} if every finite subset of
$I$ has an upper bound. To a partially ordered set $(I,\leq)$, we can associate a
category $\Zz$. The objects of $\Zz$ are the elements of $I$, and $\Hom_\Zz(i,j)=
\{a_{ji}\}$ is a singleton if $i\leq j$ and empty otherwise.\\
Let $\Aa$ be a category and $\Zz$ a category associated to a directed partially ordered set.
A functor $\ul{M}:\ \Zz\to \Aa$ will be called a {\sl direct system} with values in $\Aa$.
To $\Aa$, we associate a new category $\Aa^s$. The objects of $\Aa$ and $\Aa^s$
are the same. A morphism $M\to N$ in $\Aa^s$ is a couple $(\mu,\nu)$,
with $\mu:\ M\to N$ and $\nu:\ N\to M$ in $\Aa$ such that $\nu\circ \mu=M$,
that is, $\nu$ is a left inverse of $\mu$. A functor $\ul{M}^s:\ \Zz\to \Aa^s$
will be called a {\sl split direct system} with values in $\Aa$. We will adopt the
following notation, for all $i\leq j\in I$:
$$\ul{M}^s(i)=M_i,~~\ul{M}^s(a_{ji})=(\mu_{ji},\nu_{ij}).$$
Then $\mu_{ji}:\ M_i\to M_j$, $\nu_{ij}:\ M_j\to M_i$, and
\begin{equation}\eqlabel{3.1.1}
\nu_{ij}\circ\mu_{ji}=M_i.
\end{equation}
Consider the forgetful functor $F;\ \Aa^s\to \Aa$, $F(M)=M$, $F(\mu,\nu)=\mu$.
Then $F\circ \ul{M}^s=\ul{M}$ is a direct system with values in $\Aa$.
In \prref{3.1}, we will assume that $\colim \ul{M}=(M,\mu)$ exists. This means in particular that we have
morphisms $\mu_i:\ M_i\to M$ such that
\begin{equation}\eqlabel{3.1.2}
\mu_i=\mu_j\circ \mu_{ji}.
\end{equation}

\begin{proposition}\prlabel{3.1}
Let $\ul{M}^s:\ \Zz\to \Aa^s$ be a split direct system, and assume that
$\colim \ul{M}=(M,\mu)$ exists. Then there exist unique morphisms $\nu_i:\ M\to M_i$ in $\Aa$ such that
\begin{equation}\eqlabel{3.1.3}
\nu_i\circ \mu_i=M_i~~{\rm and}~~\nu_i=\nu_{ij}\circ \nu_j,
\end{equation}
for all $i\leq j$ in $I$.
\end{proposition}

\begin{proof}
The proof in the case where $\Aa=\Mm_A$ can be found in \cite{Joost}. In the
general case, we argue as follows.
For a fixed $i\in I$, we have a cocone $(M_i,u^i)$ on $\ul{M}$ defined as follows:
for every $k\in I$, $u^i_k:\ M_k\to M_i$ is the composition
$$\nu_{il}\circ \mu_{lk}:\ M_k\to M_l\to M_i,$$
where $l\geq i,k$. We have to show that this definition is independent of the choice
of $l$. Take $j\geq i,k$, and $m\geq l,j$. Then
$$
\nu_{im}\circ\mu_{mk}=\nu_{il}\circ \nu_{lm}\circ \mu_{ml}\circ \mu_{lk}=
 \nu_{il}\circ M_l \circ \mu_{lk}=\nu_{il}\circ \mu_{lk},
$$
and, in a similar way,
$$\nu_{im}\circ\mu_{mk}=\nu_{ij}\circ\mu_{jk}.$$
$(M_i,u^i)$ is a cocone on $\ul{M}$: take $k\geq j$ in $I$, and $l\geq i,k$; then
$$
u^i_k\circ \mu_{kj}=\nu_{il}\circ \mu_{lk}\circ \mu_{kj}=\nu_{il}\circ \mu_{lj}=u^i_j.
$$
>From the universal property of the colimit, it follows that there exists a unique
$\nu_i:\ M\to M_i$ such that 
\begin{equation}\eqlabel{3.1.4}
u^i_j=\nu_i\circ \mu_j,
\end{equation}
 for all $j\in I$. In particular,
$$\nu_i\circ \mu_i=u^i_i=\nu_{il}\circ \mu_{li}=M_i.$$
We have to show that $\nu_i=\nu_{ij}\circ \nu_j$ if $i\leq j$. To this end, it
suffices to show that
$$\nu_i\circ \mu_k=\nu_{ij}\circ \nu_j\circ \mu_k,$$
for all $k\in I$. We take $l\geq j,k$ and compute
$$\nu_i\circ \mu_k= u^i_k=\nu_{il}\circ \mu_{lk}\\
= \nu_{ij}\circ \nu_{jl}\circ\mu_{lk}=\nu_{ij}\circ u^j_k
=\nu_{ij}\circ \nu_j\circ \mu_k.
$$
We finally prove the uniqueness. Assume that $\nu'_i:\ M\to M_i$ satisfies
\equref{3.1.3}. Let $i,j\in I$, and take $k\geq i,j$. Then
$$
\nu'_i\circ\mu_j\equal{(\ref{eq:3.1.2},\ref{eq:3.1.3})}
\nu_{ik}\circ\nu'_k\circ\mu_k\circ \mu_{kj}
\equal{\equref{3.1.3}}
\nu_{ik}\circ M_k\circ\mu_{kj}=\nu_{ik}\circ \mu_{kj}=u^i_j
$$
so the $\nu'_i$ satisfy \equref{3.1.4}. By the uniqueness in the definition of colimit,
it follows that $\nu'_i=\nu_i$, for all $i\in I$.
\end{proof}

\section{Colimit comatrix corings}\selabel{4}
Let $k$ be a commutative ring. We say that a $k$-algebra $B$ has {\sl idempotent
local units} if there exists a set of idempotent elements
$\{e_i~|~i\in I\}\subset B$ such that for every finite subset $F\subset B$, there
exists $i\in I$ such that $e_ib=be_i=b$, for all $b\in F$.
If, moreover, the $e_i$ can
be chosen to be orthogonal, then we say that $B$ has {\sl orthogonal idempotent
local units}.\\
We will denote by $\Ff_k$ the category of firm $k$-algebras. 

\begin{lemma}\lelabel{4.1}
The following statements are equivalent.
\begin{enumerate}[(i)]
\item $B$ is a ring with idempotent local units;
\item There exists a split direct system $\ul{B}^s:\Zz\to\Ff_k^s$ such that $B=\colim(\ul{B},\beta)$, where $\beta_{ji}=\ul{B}(a_{ji})$ and $B_i$ is a ring with unit;
\item There exists a direct system $\ul{B}:\Zz\to\Ff_k$ such that $\colim\ul{B}=B$, where $\beta_{ji}=\ul{B}(a_{ji})$ and $B_i$ is a ring with unit.
\end{enumerate}
\end{lemma}

\begin{proof}
The statement follows immediately from Lemma 2.10 and the remark after Corollary 3.6 from \cite{Joost}. However, for sake of completeness, let us repeat a full proof using the notation we introduced in the previous section. 

\ul{$(i)\Rightarrow (ii)$}. On the index set $I$ of the idempotent local units, we define a partial ordering $\leq$ as follows: $i\leq j$ if and only if
$e_ie_j=e_je_i=e_i$. This partial ordering is directed: for all $i,j\in I$, there exists $k\in I$
such that $k\geq i,j$. Indeed, by definition of a $k$-algebra with idempotent local units, for the two elements $e_i$ and $e_j$, we can find an element $e_k$ with $k\in I$, such that $e_k$ is a local unit for both $e_i$ and $e_j$, i.e. $k\geq i,j$.
Then let $B_i=e_iBe_i$, for each $i\in I$. If $i\leq j$, then $B_i$ is a
subalgebra of $B_j$, and the inclusion map
$\beta_{ji}:\ B_i\to B_j$ is a morphism in $\Ff_k$. Also $\gamma_{ij}:\ B_j\to B_i$,
$\gamma_{ij}(b_j)=e_ib_je_i$ is a morphism of firm algebras. 
Associate to the partially ordered directed set $(I,\leq)$, a category $\Zz$ as in \seref{3}, then we have a
split direct system $\ul{B}^s:\ \Zz\to \Ff_k^s$,
$\ul{B}^s(i)=B_i$, $\ul{B}^s(a_{ji})=(\beta_{ji},\gamma_{ij})$. Clearly
$\colim \ul{B}=(B,\beta)$, with $\beta_i:\ B_i\to B$ the inclusion map.

\ul{$(ii)\Rightarrow (iii)$} is trivial.

\ul{$(iii)\Rightarrow (i)$}. Recall that module categories contain colimits and they can be described as follows.
Let
$$\Bb=\{(i,b_i)~|~i\in I,~b_i\in B_i\}.$$
be the disjoint union of the $B_i$. An equivalence relation $\sim$ on $\Bb$ is
defined as follows: $(i,b_i)\sim (j,b_j)$ if and only if there exists $k\geq i,j$
such that $\beta_{ki}(b_i)=\beta_{kj}(b_j)$, where $\beta_{ki}=B(a_{ij})$.
Then let $B=\Bb/\sim$ and $\beta_i:\ B_i\to B$, $\beta_i(b_i)=[(i,b_i)]$.
The elements of the form $[(i,1_{B_i})]$ make up a set of idempotent local units.
\end{proof}

\begin{remark}
Abrams \cite{Abrams} proved the implication (i)$\Rightarrow$(iii) under the stronger
assumption that $B$ is a ring with \emph{commuting} idempotent local units. 
\leref{4.1} tells us that the implication still holds if we drop the condition that the idempotents commute, and then we even have an equivalence. Abrams \cite[Lemma 1.5]{Abrams} 
also shows that firm modules over a ring with commuting idempotent local units
can be written as direct limits. In \leref{4.2}, this property is generalized to arbitrary
rings with idempotent local units, and it is shown that are {\emph precisely} the ones
that can be written as direct limits.
\end{remark}

Let $B$ be a $k$-algebra with idempotent local units and $A$ a $k$-algebra with unit. Let $I$ be the index set of idempotent local units of $B$ and $\Zz$ the associated category as in \leref{4.1}. 

\begin{lemma}\lelabel{4.1b}
$P$ is a firm $(B,A)$-bimodule if and only if
we can describe $P$ in the following way.
There exists a split direct system $\ul{P}^s:\ \Zz\to {\Mm_A}^s$ where we denote
for all $i\leq j\in I$:
$$\ul{P}^s(i)=P_i,~~\ul{P}^s(a_{ji})=(\sigma_{ji},\tau_{ij}),$$
and such that the following conditions hold
\begin{enumerate}[-]
\item
for all $i\leq j\in I$, $b_i\in B_i$, $p_j\in P_j$:
\begin{equation}\eqlabel{4.2.1}
\beta_{ji}(b_i)p_j=\sigma_{ji}(b_i\tau_{ij}(p_j));
\end{equation}
\item
each $P_i$ is a $(B_i,A)$-bimodule for the unital $k$-subalgebra $B_i\subset B$ and  
\item
$\colim \ul{P}=(P,\sigma)$.
\end{enumerate}
\end{lemma}

\begin{proof}
This is an immediate consequence of of \cite[Lemmas 2.7 and 2.10]{Joost}. For the sake of completeness, we give a complete proof in our present notation.

Suppose first that $P$ is a firm $(B,A)$ module. 
For each $i\in I$, we consider $P_i=e_iP$. Then
$P=\cup_{i\in I} P_i$. Moreover it is clear that $P_i$ is a left $B_i=e_iBe_i$-module and a $(B_i,A)$-bimodule.
For $i\leq j\in I$, we have right $A$-module maps
$\sigma_{ji}:\ P_i\to P_j$ (the inclusion map) and $\tau_{ij}:\ P_j\to P_i$,
$\tau_{ij}(pe_j)=pe_i$. This defines a split direct system $\ul{P}^s:\
\Zz\to \Mm_A^s$, and $\colim\ul{P}=(P,\sigma)$, with $\sigma_i:\ P_i\to P$ the
inclusion map. Finally, we check that \equref{4.2.1} holds in this situation. 
Let $i\leq j$, and take
$b_i=e_ibe_i\in B_i$, $p_j=e_jp\in P_j$. Then
$$\sigma_{ji}(b_i\tau_{ij}(p_j))=e_ibe_ie_ie_jp=e_ibe_jp=b_ie_jp=
\beta_{ji}(b_i)p_j,$$
as needed.

For the converse, the construction of the colimit is done using arguments similar to the ones in the proof of \leref{4.1}. From \prref{3.1}, we know that $\sigma_i:\ P_i\to P$ has a left inverse $\tau_i:\ P\to P_i$.
Take $p\in P$, $b\in B$. Making use of the characterisation of $B$ given in \leref{4.1}, we can find $i,j\in I$ such that
$p=\sigma_i(p_i)$, $b=\beta_j(b_j)$, with $p_i\in P_i$, $b_j\in B_j$. Take $k\geq i,j$,
and define
\begin{equation}\eqlabel{4.3.1}
bp=\sigma_k(\beta_{kj}(b_j)\sigma_{ki}(p_i)).
\end{equation}
To prove that this is a well-defined action of $B$ on $P$, we have to show that \equref{4.3.1} is independent of the choise of the index $k$. Suppose $l\ge i,j$ and consider $\sigma_l(\beta_{lj}(b_j)\sigma_{li}(p_i))$. Take any $m\ge l,k$ then we compute
\begin{eqnarray*}
\sigma_l(\beta_{lj}(b_j)\sigma_{li}(p_i))
&=&\sigma_m\circ\sigma_{ml}(\beta_{lj}(b_j)\tau_{lm}\sigma_{ml}\sigma_{li}(p_i))\\
&\stackrel{\equref{4.2.1}}{=}&\sigma_m(\beta_{ml}\beta_{lj}(b_j)\sigma_{ml}\sigma_{li}(p_i))\\
&=&\sigma_m(\beta_{mj}(b_j)\sigma_{mi}(p_i))
\end{eqnarray*}
In a similar way, we prove that we can replace $k$ by $m$ in \equref{4.3.1} and by this the left $B$ action of $P$ is independent of the choice of the index $k$.
Finally, $P$ is firm as a left $B$-module:
take $p=\sigma_i(p_i)\in P$; then $\beta_i(1_{B_i})p=p$.
\end{proof}

\begin{lemma}\lelabel{4.2}
If $P$ satisfies the equivalent conditions of \leref{4.1b}, the following formulas hold, for all $i\leq j\in I$,
$p_i\in P_i$, $p_j\in P_j$, $\varphi_i\in P_i^*$, $\varphi_j\in P_j^*$
and $b_i\in B_i$.
\begin{eqnarray}
&&(\sigma_{ji}\circ\tau_{ij})(p_j)=\beta_{ji}(1_{B_i})p_j;\eqlabel{4.2.4}\\
&&\sigma_{ji}(b_ip_i)=\beta_{ji}(b_i)\sigma_{ji}(p_i);\eqlabel{4.2.5}\\
&&\varphi_j\beta_{ji}(b_i)=((\varphi_j\circ \sigma_{ji})b_i)\circ \tau_{ij}
=\tau_{ij}^*(\sigma_{ji}^*(\varphi_j)b_i);\eqlabel{4.2.6}\\
&&\varphi_j\beta_{ji}(1_{B_i})=\varphi_j\circ\sigma_{ji}\circ\tau_{ij}=
\tau_{ij}^*(\sigma_{ji}^*(\varphi_j));\eqlabel{4.2.7}\\
&&\tau_{ij}(\beta_{ji}(b_i)p_j)=b_i\tau_{ij}(p_j);\eqlabel{4.2.8}\\
&&\tau_{ij}^*(\varphi_ib_i)=\tau_{ij}^*(\varphi_i)\beta_{ji}(b_i).\eqlabel{4.2.9}
\end{eqnarray}
\end{lemma}

\begin{proof}
\equref{4.2.4} follows after we take $b_i=1_{B_i}$ in \equref{4.2.1}. \equref{4.2.5} can
be shown as follows:
$$\beta_{ji}(b_i)\sigma_{ji}(p_i)\equal{\equref{4.2.1}}
\sigma_{ji}(b_i(\tau_{ij}\circ\sigma_{ji})(p_i))=\sigma_{ji}(b_ip_i).$$
We next prove \equref{4.2.6}. Take any $p_j\in P_j$,
\begin{eqnarray*}
&&\hspace*{-2cm}
(\varphi_j\beta_{ji}(b_i))(p_j)
= \varphi_j(\beta_{ji}(b_i)p_j)
\equal{\equref{4.2.1}} (((\varphi_j\circ\sigma_{ji})b_i)\circ\tau_{ij})(p_j)
\end{eqnarray*}
Then \equref{4.2.7} follows after we take $b_i=1_{B_i}$ in \equref{4.2.6},
and \equref{4.2.8} also follows easily:
$$\tau_{ij}(\beta_{ji}(b_i)p_j)\equal{\equref{4.2.6}}(\tau_{ij}\circ\sigma_{ji})(b_i\tau_{ij}(p_j))
=b_i\tau_{ij}(p_j).$$
\equref{4.2.9} follows immediately from \equref{4.2.8}.
\end{proof}

\begin{lemma}\lelabel{dual} 
If $P$ satisfies the equivalent conditions of \leref{4.1b}, we have a split direct system
$$\ul{P}^{*s}:\ \Zz\to {}_A\Mm^s,~~\ul{P}^{*s}(i)=P_i^*=\Hom_A(P_i,A),~~\ul{P}^{*s}(a_{ji})=
(\tau_{ij}^*,\sigma_{ji}^*),$$
where for $j\geq i$, we have defined the maps
$$\sigma_{ji}^*:\ P_j^*\to P_i^*,~~\sigma_{ji}^*(\varphi_j)=\varphi_j\circ\sigma_{ji},$$
$$\tau_{ij}^*:\ P_i^*\to P_j^*,~~\tau_{ij}^*(\varphi_i)=\varphi_i\circ\tau_{ij}.$$
Furthermore $\colim \ul{P}^{*}=(P^\dagger,\tau^\dagger)$ exists and $P^\dagger$ is a firm $(A,B)$-module.
\end{lemma}

\begin{proof}
It is straightforward to check that $P^{*s}$ is a split direct system. Since $P_i$ is a unital $(B_i,A)$-bimodule, $P_i^*$ is a unital $(A,B_i)$-bimodule. The statement follows by \leref{4.1b} using left-right duality.
\end{proof}

We will now describe the colimit of $P^{*s}$.

\begin{lemma}\lelabel{4.4}
Let $i\in I$ and $\varphi\in P^*=\Hom_A(P,A)$. There exists $\varphi_i\in P_i^*$ such that
$$\varphi=\varphi_i\circ \tau_i$$
if and only if
$$\varphi=\varphi\circ\sigma_i\circ\tau_i.$$
In this situation,
$\varphi_i$ is unique, and is given by the formula $\varphi_i=\varphi\circ\sigma_i$; furthermore,
for every $j\geq i$, $\varphi=\varphi_j\circ \tau_j$, with
$\varphi_j=\varphi_i\circ \tau_{ij}$.
\end{lemma}

\begin{proof}
If $\varphi=\varphi_i\circ\tau_i$, then $\varphi\circ\sigma_i\circ\tau_i=
\varphi_i\circ\tau_i\circ\sigma_i\circ\tau_i=\varphi_i\circ\tau_i=\varphi$. The converse
is obvious. If $\varphi=\psi\circ\tau_i$, then $\varphi\circ\sigma_i=
\psi\circ\tau_i\circ\sigma_i=\psi$. If $j\geq i$, then
$\varphi_i\circ \tau_{ij}\circ\tau_j=\varphi_i\circ\tau_i=\varphi$.
\end{proof}

Let $P^\dagger=\{\varphi\in P^*~|~\exists i\in I:~\varphi=\varphi\circ\sigma_i\circ\tau_i\}$. More explicit, using the characterisations \leref{4.1} and \leref{4.1b} we get $P^\dagger=\{\varphi\in P^*~|~\exists i\in I:~\varphi(p)=\varphi(e_ip), {\rm ~for~all~}p\in P\}$.
For every $i\in I$, we have a map 
$$\tau_i^*:\ P_i^*\to P^\dagger,~~\tau_i^*(\varphi_i)=\varphi_i\circ\tau_i.$$

\begin{proposition}\prlabel{4.5}
With notation as above, $\colim \ul{P}^*=(P^\dagger, \tau^*)$.
\end{proposition}

\begin{proof}
First, $(P^\dagger, \tau^*)$ is a cocone on $Q^*$
since, for all $i\leq j$ and $\varphi_i\in P_i^*$, we have
$$(\tau_j^*\circ \tau_{ij}^*)(\varphi_i)=\varphi_i\circ \tau_{ij}\circ \tau_j
=\varphi\circ\tau_i=\tau_i^*.$$
Let $(M,m)$ be another cocone on $Q^*$. This means that $m_i:\ P_i^*\to M$
and $m_j\circ\tau_{ij}^*=m_i$ if $i\leq j$. We then define $f:\ P^\dagger\to M$
as follows: $f(\varphi_i\circ\tau_i)=m_i(\varphi_i)$, for every $i\in I$ and $\varphi_i\in P_i^*$.
Let us show that $f$ is well-defined. Assume that
$$\varphi=\varphi_i\circ\tau_i=\varphi_j\circ\tau_j.$$
Take $k\geq i,j$. Then $\varphi=\varphi_k\circ\tau_k$ with $\varphi_k=\varphi_i\circ \tau_{ik}$
(see \leref{4.4}). Then
$$m_k(\varphi_k)=m_k(\varphi_i\circ \tau_{ik})=(m_k\circ \tau_{ik}^*)(\varphi_i)=
m_i(\varphi_i).$$
In a similar way, we have that $m_j(\varphi_j)=m_k(\varphi_k)$, and it follows that $f$
is well-defined.
Finally, $(f\circ \tau_i^*)(\varphi_i)=f(\varphi_i\circ \tau_i)=m_i(\varphi_i)$.
\end{proof}

The right $B$-action on $P^\dagger$ can be described as follows: take $\varphi=
\varphi_i\circ\tau_i\in P^\dagger$ and $b=\beta_j(b_j)\in B$. For $k\geq i,j$, we have
\begin{equation}\eqlabel{4.5.1}
\varphi b=((\varphi_i\circ\tau_{ik})\beta_{kj}(b_j))\circ\tau_k.
\end{equation}
In particular, we have, for $\varphi_i\in P_i^*$ and $b_i\in B_i$:
\begin{equation}\eqlabel{4.5.2}
(\varphi_i\circ\tau_i)\beta_i(b_i)=(\varphi_i b_i)\circ\tau_i.
\end{equation}
In explicit form this means $(\varphi b)(p)=\varphi_i(e_ie_jbe_jp)$ or just $(\varphi b)(p)=\varphi(bp)$.

\begin{lemma}\lelabel{4.5a}
If $P$ satisfies the equivalent conditions of \leref{4.1b}, then we have 
for all $i\in I$, $b_i\in B_i$, $p\in P$ and $\varphi\in P^\dagger$,
\begin{eqnarray}
\beta_i(b_i)p&=&\sigma_i(b_i\tau_i(p));\eqlabel{4.5a.1}\\
\varphi\beta_i(b_i)&=& (\varphi\circ\sigma_i)b_i\circ\tau_i.\eqlabel{4.5a.2}
\end{eqnarray}
\end{lemma}

\begin{proof}
By the characterisation of \leref{4.1} and \leref{4.1b}, we can write $b_i=e_ibe_i$ and $\tau_i(p)=e_ip$, where $e_i$ is an idempotent in $B$. Moreover the maps $\beta_i$ and $\sigma_i$ are injections. With this information in hand we easily find
\begin{eqnarray*}
\beta_i(b_i)p&=&e_ibe_ip\\
=\sigma_i(b_i\tau_i(p))&=&e_ibe_ie_ip.
\end{eqnarray*}
The other equation follows by 
$$\varphi\beta_i(b_i)(p)=\varphi(\beta_i(b_i)p)
=\varphi(\sigma_i(b_i\tau_i(p)))=(\varphi\circ\sigma_i)b_i\circ\tau_i(p),$$
where we used \equref{4.5a.1} in the second equality.
\end{proof}

\begin{proposition}\prlabel{4.6}
We have a directed system $\ul{G}:\ \Zz\to {}_A\Mm_A$, $G(i)=P_i^*\ot_{B_i}P_i$, and
$$\ul{G}(a_{ji}):\ P_i^*\ot_{B_i}P_i\to P_j^*\ot_{B_j}P_j,~~
\ul{G}(a_{ji})(\varphi_i\ot_{B_i}p_i)=
\varphi_i\circ \tau_{ij}\ot_{B_j}\sigma_{ji}(p_i).$$
\end{proposition}

\begin{proof}
We first show that $\ul{G}(a_{ji})$ is well-defined. For all $\varphi_i\in P_i^*$, $p_i\in P_i$
and $b_i\in B_i$, we have
\begin{eqnarray*}
&&\hspace*{-2cm}\ul{G}(a_{ji})(\varphi_i\ot_{B_i}b_i\cdot p_i)=
\varphi_i\circ \tau_{ij}\ot_{B_j}\sigma_{ji}(b_ip_i)\\
&\equal{\equref{4.2.5}}& \tau_{ij}^*(\varphi_i)\ot_{B_j} \beta_{ji}(b_i)\sigma_{ji}(p_i)
= \tau_{ij}^*(\varphi_i)\beta_{ji}(b_i)\ot_{B_j} \sigma_{ji}(p_i)\\
&\equal{\equref{4.2.9}}&
\tau_{ij}^*(\varphi_ib_i)\ot_{B_j} \sigma_{ji}(p_i)
= \ul{G}(a_{ji})(\varphi_ib_i\ot_{B_i}p_i).
\end{eqnarray*}
If $i\leq j\leq k$, then we have
\begin{eqnarray*}
&&\hspace*{-2cm}(\ul{G}(a_{kj})\circ \ul{G}(a_{ji}))(\varphi_i\ot_{B_i}t_i\cdot p_i)=
\varphi_i\circ\tau_{ji}\circ\tau_{jk}\ot_{B_k} (\sigma_{kj}\circ\sigma_{ji})(p_i)\\
&=&\varphi_i\circ \tau_{ik}\ot_{B_k}\sigma_{ki}(p_i)
=\ul{G}(a_{ki})(\varphi_i\ot_{B_i}p_i).
\end{eqnarray*}
\end{proof}

Let $P$ be a module satisfying the equivalent conditions of \leref{4.1b}. Suppose that $P_i$ is finitely generated and projective as right $A$-module for all $i\in I$. Let $E_i=\sum z_i\ot_A z_i^*$ be a finite dual basis of $P_i\in \Mm_A$;
we omitted the summation index. $E_i$ is the unique element of $P_i\ot_AP_i^*$
satisfying the formulas
\begin{equation}\eqlabel{4.2.3}
p_i=\sum z_iz_i^*(p_i)~~;~~\varphi_i=\sum z_i^*\varphi(z_i),
\end{equation}
for all $p_i\in P_i$ and $\varphi_i\in P_i^*$.
With these notation, we have the following lemma.

\begin{lemma}\lelabel{4.7}
\begin{enumerate}[(i)]
\item For all $b_i\in B_i$,
\begin{equation}\eqlabel{4.11.1}
\sum b_iz_i\ot_Az_i^*=\sum z_i\ot_Az_i^*b_i,
\end{equation}
\item If $i\leq j$, then
\begin{equation}\eqlabel{4.7.1}
E_i=\sum \tau_{ij}(z_j)\ot_A z_j^*\circ \sigma_{ji}=e_iz_j\ot_A z_{j|P_i}.
\end{equation}
\end{enumerate}
\end{lemma}

\begin{proof}
\ul{$(i)$}. This follows from the fact that $P_i$ is a 
$(B_i,A)$-bimodule.

\ul{$(ii)$}.
We show that the right hand side of \equref{4.7.1} satisfies \equref{4.2.3}. For all
$p_i\in P_i$, we have
$$
\sum \tau_{ij}(z_j)(z_j^*\circ \sigma_{ji})(p_i)=
\sum \tau_{ij}\Bigl(z_j z_j^*\bigl(\sigma_{ji}(p_i)\bigr)\Bigr)=
\tau_{ij}(\sigma_{ji}(p_i))=p_i.$$
\end{proof}

For the remaining part of this paper, we will concentrate on modules that are locally projective in the sense of \'Anh and M\'arki \cite{AM} (strongly locally projective in the terminology of \cite{Joost}). We will need a more restrictive characterisation than \leref{4.1b}. Recall first the definition of a morphism $\eta:B\to B'$ of rings with (idempotent) local units. This is a ringmorphism $\eta$ satisfying the property that for every finite subset $F'\subset B'$, we can find an (idempotent) local unit $e_i\in B$ such that $\eta(e_i)$ is an (idempotent) local unit for all elements of $F'$.

\begin{lemma}\lelabel{4.11}
The following statements are equivalent
\begin{enumerate}[(i)]
\item
$P$ satisfies the equivalent condtions of \leref{4.1b}, in addition $P_i$ is finitely generated and projective as right $A$-module for all $i\in I$ and $\colim \ul{P^*}=(P^\dagger,\tau^\dagger)$.
\item $S=P\otimes_A P^\dagger$ is a ring with idempotent local units, $P$ is a firm left $S$-module, $P^\dagger$ is a firm right $S$-module and there exists a unique morphism of rings with idempotent local units
$\eta:\ B\to
P\ot_AP^\dagger$ such that
\begin{equation}\eqlabel{4.11.2}
\eta(\beta_i(b_i))=\sum \sigma_i(b_iz_i)\ot_A z_i^*\circ \tau_i=
\sum \sigma_i(z_i)\ot_A z_i^*b_i\circ \tau_i,
\end{equation}
for all $i\in I$, $b_i\in B_i$.
\item
$P$ is strongly $P^\dagger$-locally projective as right $A$-module and $P^\dagger$ is strongly $P$-locally projective as left $A$-module. $P$ is a firm left $B$-module and $P^\dagger$ is a firm right $B$-module.
\end{enumerate}
\end{lemma}

\begin{proof}
\ul{$(i)\Rightarrow (ii)$}
By \leref{dual} $\ul{P^{*s}}$ is a split direct system and obviously $P^*_i$ is finitely generated an projective as left $A$-module for every $i\in I$.
The first part of statement $(ii)$ follows now from \cite[Corollary 3.6]{Joost}.
The second equality in \equref{4.11.2} is an immediate consequence of \equref{4.11.1}.
Let us show that $\eta$ is well-defined. Take $b\in B$, and assume that
$b=\beta_i(b_i)=\beta_j(b_j)$, for some $i,j\in I$, $b_i\in B_i$, $b_j\in B_j$.
Take $k\geq i,j$, and let $b_k=\beta_{ki}(b_i)=\beta_{kj}(b_j)$. We compute
\begin{eqnarray*}
&&\hspace*{-2cm}
\sum \sigma_i(b_iz_i)\ot_A z_i^*\circ \tau_i\equal{\equref{4.7.1}}
\sum (\sigma_k\circ\sigma_{ki})(b_i\tau_{ik}(z_k))\ot_A
z_k^*\circ\sigma_{ki}\circ\tau_{ik}\circ\tau_k\\
&\equal{(\ref{eq:4.2.1},\ref{eq:4.2.8})}&
\sum \sigma_k(\beta_{ki}(b_i)z_k)\ot_A (z_k^*\beta_{ki}(1_{B_i}))\circ \tau_k\\
&\equal{\equref{4.11.1}}&
\sum \sigma_k(\beta_{ki}(1_{B_i})\beta_{ki}(b_i)z_k)\ot_A z_k^*\circ \tau_k\\
&=& \sum \sigma_k(\beta_{ki}(1_{B_i}b_i)z_k)\ot_A z_k^*\circ \tau_k\\
&=& \sum \sigma_k(b_kz_k)\ot_A z_k^*\circ \tau_k.
\end{eqnarray*}
In a similar way, we prove that
$$\sum \sigma_j(b_jz_j)\ot_A z_j^*\circ \tau_j=\sum \sigma_k(b_kz_k)\ot_A z_k^*\circ \tau_k,$$
and it follows that the right hand side of \equref{4.11.2} is independent of the choice of $i$.
Next we prove that $\eta$ is a ringmorphism.  Take two elements $b, b'\in B$ and choose $i$ big enough such that $b=\beta_i(b_i)$ and $b'=\beta_i(b'_i)$. Let us denote $E_i=\sum z_i\ot_Az_i^*=\sum \tilde{z}_i\ot_A\tilde{z}_i^*$
\begin{eqnarray*}
\eta(b)\eta(b')&=&\sum \sigma_i(b_iz_i) z_i^*\circ \tau_i\circ\sigma_i(b'_i\tilde{z}_i)\ot_A \tilde{z}_i^*\circ \tau_i\\
&=&\sum \sigma_i(b_i b'_i\tilde{z}_i)\ot_A\tilde{z}_i^*\circ\tau=\eta(bb')
\end{eqnarray*}
Finally, the idempotent local units of $P\ot_A P^\dagger$ are of the form $\sum \sigma_i(z_i)\ot_A z_i^*\circ \tau_i$, these are exactly given by $\eta(1_{B_i})$, so $\eta$ is a morphism of rings with idempotent local units.

\ul{$(ii)\Rightarrow (iii)$}. By \cite[Corollary 3.6]{Joost} we only have to prove that $P$ and $P^\dagger$ are firm $B$-modules under the action induced by the morphism $\eta$. This is a consequence of the fact that $\eta$ is a morphism of rings with enough idempotents. Take $p\in P$, then we know there exists an idempotent $e\in B$ such that $\eta(e)\in P\otimes_AP^\dagger$ is a local unit for $p$. Thus $e \cdot p= \eta(e) p=p$ and $P$ is a firm $B$-module. Analougously one proves $P^\dagger$ is a firm right $B$-module. 

\ul{$(iii)\Rightarrow (i)$}. Follows from \cite[Corollary 3.6]{Joost} and \leref{4.1b}.
\end{proof}

For every $i\in I$, consider bimodule maps
\[
\begin{array}{rcl}
\coev_{P_i}:& B_i\to P_i\otimes_AP_i^*,&\quad \coev_{P_i}(b_i)=b_iE_i=E_ib_i\\
\ev_{P_i}:&P_i^*\otimes_{B_i}P_i,&\quad\ev_{P_i}(\varphi_i\otimes_{B_i}p_i)=\varphi_i(p_i)
\end{array}
\]
Then $(B_i,A,P_i,P_i^*,\coev_{P_i},\ev_{P_i})$ is a comatrix coring
context, so we have a comatrix coring $(\ul{G}(i),\Delta_i,\varepsilon_i)$ with
$$\Delta_i(\varphi_i\ot_{B_i} p_i)=\varphi_i\ot_{B_i}E_i\ot_{B_i} p_i~~{\rm and}~~
\varepsilon_i(\varphi_i\ot_{B_i} p_i)=\varphi_i(p_i).$$
$\ul{G}(i)$ is a finite comatrix coring, as introduced in \cite{Kaoutit}.

\begin{proposition}\prlabel{4.8}
Suppose the equivalent conditions of \leref{4.11} hold, and consider the directed system $\ul{G}$ from \prref{4.6}.
Then $(\ul{G},\Delta,\varepsilon)$ is a coalgebra in ${\rm Func}(\Zz,{}_A\Mm_A)$.
\end{proposition}

\begin{proof}
It suffices to show that $\Delta$ and $\varepsilon$ are natural transformations, or,
equivalently, that $\ul{G}(a_{ji})$ is a morphism of corings, for every $i\leq j$, or
$$(\ul{G}(a_{ji})\ot_A \ul{G}(a_{ji}))\circ\Delta_i=\Delta_j\circ \ul{G}(a_{ji})~~;~~
\varepsilon_i=\varepsilon_j\circ \ul{G}(a_{ji}).$$
For all $\varphi_i\in P_i^*$ and $p_i\in P_i$, we compute
\begin{eqnarray*}
&&\hspace*{-2cm}
(\Delta_j\circ \ul{G}(a_{ji}))(\varphi_i\ot_{B_i}p_i)=
\Delta_j\bigl(\varphi_i\circ \tau_{ij}\ot_{B_j}\sigma_{ji}(p_i)\bigr)\\
&=& \varphi_i\circ \tau_{ij}\ot_{B_j}E_j\ot_{B_j}\sigma_{ji}(p_i)\\
&=& \sum \varphi_i\circ \tau_{ij}\circ \sigma_{ji}\circ \tau_{ij}
\ot_{B_j}z_j\ot_Az_j^*\ot_{B_j}
(\sigma_{ji}\circ \tau_{ij}\circ \sigma_{ji})(p_i)\\ 
&\equal{(\ref{eq:4.2.4},\ref{eq:4.2.8})}& 
\sum \varphi_i\circ \tau_{ij}
\ot_{B_j}(\sigma_{ji}\circ \tau_{ij})(z_j)\ot_Az_j^*\circ \sigma_{ji}\circ \tau_{ij}\ot_{B_j}
\sigma_{ji}(p_i)\\
&\equal{\equref{4.7.1}}&
\sum \varphi_i\circ \tau_{ij}
\ot_{B_j}\sigma_{ji}(z_i)\ot_Az_i^*\circ \tau_{ij}\ot_{B_j}
\sigma_{ji}(p_i)\\
&=& \sum  \ul{G}(a_{ji})(\varphi_i\ot_{B_i}z_i)\ot_A  \ul{G}(a_{ji})(z_i^*\ot_{B_i}p_i)\\
&=&\bigl(( \ul{G}(a_{ji})\ot_A  \ul{G}(a_{ji}))\circ\Delta_i\bigr)(\varphi_i\ot_{B_i}p_i)
\end{eqnarray*}
and
$$
\varepsilon_j\bigl( \ul{G}(a_{ji})(\varphi_i\ot_{B_i}p_i)\bigr)=
(\varphi_i\circ\tau_{ij}\circ\sigma_{ji})(p_i)
= \varphi_i(p_i)=\varepsilon_i(\varphi_i\ot_{B_i}p_i).$$
\end{proof}

\begin{proposition}\prlabel{4.9}
Under the same conditions as \prref{4.8}, $\colim \ul{G}=(P^\dagger\ot_{B}P,g)$, with
$$g_i:\ \ul{G}(i)=P_i^*\ot_{B_i}P_i\to P^\dagger\ot_{B}P,~~
g_i(\varphi_i\ot_{B_i}p_i)=\varphi_i\circ\tau_i\ot_{B}\sigma_i(p_i).$$
\end{proposition}

\begin{proof}
We first show that $g_i$ is well-defined. For all $b_i\in B_i$, we have
\begin{eqnarray*}
&&\hspace*{-15mm}
g_i(\varphi_i b_i\ot p_i)
=\varphi_i b_i\circ\tau_i \ot_B \sigma_i( p_i)
= (\varphi_i\circ\tau_i\circ\sigma_i) b_i\circ\tau_i \ot_B \sigma_i( p_i)\\
&\equal{\equref{4.5a.2}}&
(\varphi_i\circ\tau_i)\beta_i(b_i)\ot_B \sigma_i( p_i)
= (\varphi_i\circ\tau_i)\ot_B \beta_i(b_i)\sigma_i( p_i)\\
&\equal{\equref{4.5a.1}}&
(\varphi_i\circ\tau_i)\ot_B \sigma_i(b_i\tau_i(\sigma_i p_i))
=(\varphi_i\circ\tau_i)\ot_B \sigma_i(bp_i)
= g_i(\varphi_i \ot b_ip_i).
\end{eqnarray*}
Let us now prove that $(P^\dagger\ot_{B}P,g)$ is a cocone on $\ul{G}$.
Indeed, if $i\leq j$, then
\begin{eqnarray*}
&&\hspace*{-15mm}
(g_j\circ \ul{G}(a_{ji}))(\varphi_i\ot_{B_i}p_i)=
g_j(\varphi_i\circ\tau_{ij}\ot_{B_j}\sigma_{ji}(p_i))\\
&=&\varphi_i\circ\tau_{ij}\circ \tau_j \ot_{B}\sigma_j(\sigma_{ji}(p_i))
= \varphi_i\circ\tau_i \ot_{B}\sigma_i(p_i)=g_i(\varphi_i\ot_{B_i}p_i).
\end{eqnarray*}
Let $(M,m)$ be another cocone on $\ul{G}$. Then $m_i:\ P_i^*\ot_{B_i}P_i\to M$ and
$m_j\circ \ul{G}(a_{ji})=m_i$ if $j\geq i$. 
We define $f:\ P^\dagger\ot P
\to M$ as follows. For $\varphi\in P^\dagger$ and $p\in P$, we can find
$i\in I$, $\varphi_i\in P_i^*$ and $p_i\in P_i$ such that $p=\sigma_i(p_i)$
and $\varphi=\varphi_i\circ\tau_i$; we then define
$$f(\varphi\ot p)=m_i(\varphi_i\ot_{T_i}p_i).$$
We have to show that $f$ is well-defined. If $k\geq i$, then we have that
$\varphi=\varphi_k\circ \tau_k$ and $p=\sigma_k(p_k)$ with
$\varphi_k=\varphi_i\circ\tau_{ik}$ and $p_k=\sigma_{ki}(p_i)$. We then find that
\begin{eqnarray*}
&&\hspace*{-2cm}
m_k(\varphi_k\ot_{B_k}p_k)=m_k(\varphi_i\circ\tau_{ik}\ot_{B_k}\sigma_{ki}(p_i))\\
&=& (m_k\circ \ul{G}(a_{ki}))(\varphi_i\ot_{B_i}p_i)=m_i(\varphi_i\ot_{B_i}p_i).
\end{eqnarray*}
We will now show that $f$ induces a map $f:\ P^\dagger\ot_{B} P
\to M$. To this end, we need to prove that
$$f(\varphi b\ot p)=f(\varphi\ot b p),$$
for all $\varphi\in P^\dagger$, $p\in P$ and $b\in B$. We can find
$i\in I$, $b_i\in B_i$, $\varphi\in P_i^*$ and $p_i\in P_i$ such that 
$b=\beta_i(b_i)$, $p=\sigma_i(p_i)$
and $\varphi=\varphi_i\circ\tau_i$. Then we compute that
\begin{eqnarray*}
&&\hspace*{-15mm}
f(\varphi b\ot p)=f((\varphi_i\circ\tau_i)\beta_i(b_i)\ot\sigma_i(p_i))\\
&\equal{\equref{4.5a.2}}&
f(((\varphi_i\circ\tau_i\circ\sigma_i)b_i)\circ\tau_i\ot \sigma_i(p_i))
=f(\varphi_ib_i\circ\tau_i\ot \sigma_i(p_i))\\
&=& m_i(\varphi_i b_i\ot_{B_i} p_i)
= m_i(\varphi_i\ot_{B_i} b_ip_i)\\
&=& f(\varphi_i\circ\tau_i \ot \sigma_i(b_ip_i))
= f(\varphi_i\circ\tau_i \ot \sigma_i(b_i\tau_i(\sigma_i(p_i))))\\
&\equal{\equref{4.5a.1}}&
f(\varphi_i\circ\tau_i \ot \beta_i(b_i)\sigma_i(\tau_i(\sigma_i(p_i))))
= f(\varphi\ot b p).
\end{eqnarray*}
Finally,
$$f(g_i(\varphi_i\ot_{B_i}p_i))=f(\varphi_i\circ\tau_i\ot_{B}\sigma_i(p_i))
=m_i(\varphi_i\ot_{B_i}p_i).$$
\end{proof}

The following result now follows immediately from Propositions \ref{pr:2.1},
 \ref{pr:4.8} and  \ref{pr:4.9}.

\begin{corollary}\colabel{4.10}
If the equivalent conditions of \leref{4.11} hold,
$\Gg=P^{\dagger}\ot_{B}P$ is an $A$-coring, with comultiplication
and counit given by the following formulas, for all $i\in I$, $\varphi_i\in P_i^*$
and $p_i\in P_i$:
$$\Delta(\varphi_i\circ \tau_i\ot_{B}\sigma_i(p_i))=
\sum \Delta(\varphi_i\circ \tau_i\ot_{B}\sigma_i(z_i)\ot_A
z_i^*\circ \tau_i\ot_{B}\sigma_i(p_i)),$$
$$\varepsilon(\varphi_i\circ \tau_i\ot_{B}\sigma_i(p_i))=\varphi_i(p_i).$$
As before, $E_i=\sum z_i\ot_A z_i^*$ is the finite dual basis of
$P_i\in \Mm_A$.
\end{corollary}

We will now show that $P^{\dagger}\ot_{B}P$ can be constructed starting from
a G\' omez-Vercruysse comatrix coring context, as described in \seref{1}.
We already know that $P$ and $P^\dagger$ are firm bimodules.

\begin{proposition}\prlabel{4.12}
If the equivalent conditions of \leref{4.11} hold,
$(B,A,P,P^\dagger,\eta,\varepsilon)$ is a comatrix coring
context, where $\varepsilon:\ P^{\dagger}\ot_BP\to A$ is the restriction of the evaluation map $P^*\ot_BP\to A$.
\end{proposition}

\begin{proof}
We have to show that \equref{1.1.1} holds. Take $b=\beta_i(b_i)\in B$,
$p=\sigma_i(p_i)\in P$ and $\varphi=\varphi_i\circ\tau_i\in P^\dagger$. Then
\begin{eqnarray*}
&&\hspace*{-2cm}
b^{-}\varepsilon(b^+\ot_A p)=
\sum \sigma_i(b_iz_i)(z_i^*\circ\tau_i)(\sigma_i(p_i))\\
&=&\sum \sigma_i(b_iz_iz_i^*(p_i))=\sigma_i(b_ip_i)
=\beta_i(b_i)\sigma_i(p_i)=bp,
\end{eqnarray*}
and
\begin{eqnarray*}
&&\hspace*{-15mm}
\bigl(\varepsilon(\varphi\ot_A b^-)b^+\bigr)(p)=
\bigl(\sum \varphi(\sigma_i(z_i))(z_i^*b_i)\circ\tau_i\bigr)(p)
= \sum \varphi(\sigma_i(z_i))z_i^*(b_i\tau_i(p))\\
&=& \sum \varphi(\sigma_i(z_iz_i^*(b_ip_i)))
= \varphi(\sigma_i(b_ip_i))=\varphi(bp)=(\varphi b)(p),
\end{eqnarray*}
hence
$\varepsilon(\varphi\ot_A b^-)b^+=\varphi b$.
\end{proof}

\begin{example}
Let $B$ be a $k$-algebra with {\sl orthogonal}
idempotent local units and let $\{e_i~|~i\in I\}$ be a complete set of idempotents. 
For all $i,j\in I$, let $B_{ij}=e_iBe_j$. Then
$B=\bigoplus_{i,j\in I} B_{ij}$,
and a firm left $B$-module $P$ can then be written as
$P=\bigoplus_{i\in I} P_i$,
with $P_i=e_iP$ a left $B_i=B_{ii}$-module. For each ${i\in I}$, we take a $(B_i,A)$-bimodule
$P_i$ which is finitely generated and projective as a right $A$-module, and
we put $P=\bigoplus_{i\in I} P_i$. It is not hard to see that
$P^\dagger=\bigoplus_{i\in I} P_i^*$,
and we have a comatrix coring
$P^\dagger\ot_B P$
. This way we recover the comatrix corings that were considered first in \cite[Proposition 5.2]{KGT}.
\end{example}

\begin{example}\exlabel{5.1}
As a special case of the previous example, consider now the case where the orthogonal idempotents are central in $B$, then the situation simplifies to $B=\oplus_{i\in I}B_i$, where $B_i=Be_i$.
\end{example}

The functor $K:\ \Mm_B\to \Mm^\Gg$
can be described as follows. Take $M\in \Mm_B$, and, as in \leref{4.1b}, let 
$M_i=Me_i$. We have a split direct system $\ul{F}:\ \Zz\to \Mm_A^s$:
$$\ul{F}(i)=M_i\ot_{B_i}P_i~~;~~
\ul{F}(a_{ji})=(\mu_{ji}\ot \sigma_{ji},\nu_{ij}\ot\tau_{ij}).$$
Then
$K(M)=\colim \ul{F}$, with the obvious coaction.

In view of \thref{1.3}, it is important to know when $P\in {}_B\Mm$ is
(faithfully) flat. We have the following results.

\begin{proposition}\prlabel{5.2}
Let $B$ be a $k$-algebra with idempotent local units, and take $P\in {}_B\Mm$.
If for every $i\in I$, there exists $j\geq i$ such that $P_j\in {}_{B_j}\Mm$ is flat, then
$P\in {}_B\Mm$ is flat.
\end{proposition}

\begin{proof}
Let $f:\ N'\to N$  be an injective map in $\Mm_B$, and $x\in \ker(f\ot_B P)$.
$N'\ot_B P$ is the colimit of the $N'_i\ot_{B_i} P_i$, so $x$ can be represented
by $\sum_r n'_r\ot_{B_i} p_r$ with $n'_r\in N'_i$, $p_r\in P_i$.
$\sum_r f(n'_r)\ot_{B_i} p_r$ represents zero in $N\ot_B P$, so, replacing $i$ by
a bigger index, we can assume that $\sum_r f(n'_r)\ot_{B_i} p_r=0\in N_i\ot_{B_i} P_i$.
Replace $i$ by a bigger index such that $P_i\in {}_{B_i}\Mm$ is flat.
Then $\sum_r n'_r\ot_{B_i} p_r=0$ in $N'_i\ot_{B_i} P_i$, and this implies that $x=0$.
\end{proof}

\begin{proposition}\prlabel{5.3}
Let $B$ be a $k$-algebra with idempotent local units, and assume that
$P\in {}_B\Mm$ is (faithfully) flat. If $i\in I$ is such that $e_i$ is central in $B$, then
$P_i$ is (faithfully) flat as a left $B_i$-module.
\end{proposition}

\begin{proof}
Take $N\in \Mm_{B_i}$. We have $\gamma_i:\ B\to B_i$, making $N\in \Mm_B$
via restriction of scalars. Then we claim that we have an isomorphism of
$k$-modules
\begin{equation}\eqlabel{5.2.1}
N\ot_BP\cong N\ot_{B_i}P_i.
\end{equation}
Indeed, the map
$$f:\ N\ot_{B_i}P_i\to N\ot_BP,~~f(n\ot_{B_i}p_i)=n\ot_B p_i$$
has an inverse $g$ given by
$$g(n\ot_B p)=n\ot_{B_i} e_i p.$$
$g$ is well-defined since
$$g(nb\ot p)=g(ne_ib\ot p)=g(ne_ibe_i\ot p)=ne_ibe_i\ot_{B_i} e_ip=
n\ot_{B_i} e_ibe_ip=g(n\ot bp).$$
Assume that $P\in {}_B\Mm$ is faithfully flat. A sequence
$$0\to N'\to N\to N''\to 0$$
is exact in $\Mm_{B_i}$ if and only if
$$0\to N'\ot_BP\to N\ot_BP\to N''\ot_BP\to 0$$
is exact in $\Mm_k$, and, by \equref{5.2.1}, this is equivalent to
exactness of the sequence
$$0\to N'\ot_{B_i}P_i\to N\ot_{B_i}P_i\to N''\ot_{B_i}P_i\to 0.$$
\end{proof}

We remark that the condition that $e_i$ is central is fulfilled in the situation of
\exref{5.1}.
The condition that the $e_i$ are central is also needed in the proof
of our next result. We have seen that the comatrix coring is the colimit of the directed system
$\ul{G}$ discussed in \prref{4.6}. If we work over an algebra with central idempotent
local units, then this system is split.

\begin{proposition}\prlabel{5.4}
Let $B$ be a $k$-algebra with central idempotent local units and 
suppose the equivalent conditions of \leref{4.11} hold, then the direct system $\ul{G}$ of \prref{4.6} splits.
$\ul{G}^s:\ \Zz\to {}_A\Mm_A^s$, with
$\ul{G}^s(a_{ji})=(g_{ji},h_{ij})$,
where 
$$h_{ij}(\varphi_j\ot_{B_j}p_j)=\varphi_j\circ\sigma_{ji}\ot_{B_i}\tau_{ij}(p_j)=\varphi_{j|P_i}\ot_{B_i}e_ip_j,$$
for all $\varphi_j\in P_j^*$ and $p_j\in P_j$.
\end{proposition}

\begin{proof}
Let us show that $h_{ij}$ is well-defined; all the rest is obvious. First we compute
for $\varphi_j\in P_j^*$, $b_j\in B_j$ and $p_i\in P_i$
that
$$
(\varphi_j b_j)(\sigma_{ji}p_i)=\varphi_j(b_jp_i)=\varphi_j(b_je_ip_i)
= \varphi_j(e_ib_je_ip_i)=(\varphi_j\circ \sigma_{ji})(e_ib_je_i)(p_i),$$
where we used the fact that $e_i$ is central. Then we compute
\begin{eqnarray*}
&&\hspace*{-15mm}
h_{ij}(\varphi_jb_j\ot p_j)=\varphi_jb_j\circ\sigma_{ji}\ot_{B_i} e_ip_j=
(\varphi_j\circ \sigma_{ji})(e_ib_je_i)\ot_{B_i} e_ip_j\\
&=& \varphi_j\circ \sigma_{ji}\ot_{B_i}e_ib_je_ie_ip_j=
\varphi_j\circ \sigma_{ji}\ot_{B_i}e_ib_jp_j=h_{ij}(\varphi_j\ot b_jp_j).
\end{eqnarray*}
\end{proof}

\section{Factorizing split direct systems}\selabel{6}
In this Section, we consider split direct systems $\ul{P}^s:\ \Zz\to \Mm_{A, \rm fgp}^s$
that factorize through a $k$-linear category $\Aa$: we assume that there
exists a split direct system
$$\ul{M}^s:\ \Zz\to \Aa^s,~~,\ul{M}^s(i)=M_i,~~
\ul{M}^s(a_{ji})=(\mu_{ji},\nu_{ij})$$
and a functor $\omega:\ \Aa\to \Mm_A$ such that $\ul{P}^s=\omega\circ\ul{M}^s$,
or
$$P_i=\omega(M_i),~~\sigma_{ji}=\omega(\mu_{ji}),~~\tau_{ij}=\omega(\nu_{ij}).$$
For every $i\in I$,
$T_i=\End_\Aa(M_i)$
is a $k$-algebra with unit. For $i\leq j$, we have a multiplicative map
$$\rho_{ji}:\ T_i\to T_j,~~\rho_{ji}(t_i)=\mu_{ji}\circ t_i\circ\nu_{ij}.$$
This defines a direct system $\ul{T}:\ \Zz\to \Ff_k$, $\ul{T}_i=T_i$,
$\ul{T}(a_{ji})=\rho_{ji}$. If $t_i\in T_i=\End_\Aa(M_i)$, then $\omega(t_i)
\in \End_A(P_i)$. Hence
$P_i$ is a $(T_i,A)$-bimodule, with left $T_i$-action
given by
$$t_i\cdot p_i=\omega(t_i)(p_i).$$
We claim that \equref{4.2.1} holds. Indeed, for all $i\leq j$, $t_i\in T_i$ and
$p_j\in P_j$, we have
\begin{eqnarray*}
&&\hspace*{-2cm}
\rho_{ji}(t_i)\cdot p_j=(\mu_{ji}\circ t_i\circ \nu_{ij})\cdot p_j
=\omega(\mu_{ji}\circ t_i\circ \nu_{ij})(p_j)\\
&=& (\sigma_{ji}\circ \omega(t_i)\circ \tau_{ij})(p_j)
=\sigma_{ji}(t_i\cdot \tau_{ij}(p_j)).
\end{eqnarray*}
Applying the results of \seref{4}, we obtain a comatrix coring. We will now assume
that $\colim\ul{M}=(M,\mu)$ exists, and that $\omega$ preserves colimits.
We will give an explicit description of $\colim \ul{T}$, and provide some alternative
descriptions of the comatrix coring. Using \prref{3.1} , we obtain morphisms
$\nu_i:\ M\to M_i$. Let $\sigma_i=\omega(\mu_i)$, $\tau_i=\omega(\nu_i)$.
We consider the $k$-algebra
$T=\End_\Aa(M)$. For every $i\in I$, $e_i=\mu_i\circ \nu_i$
is an idempotent in $T$. We also have
\begin{equation}\eqlabel{6.1.1}
e_i\circ\mu_i=\mu_i~~{\rm and}~~\nu_i\circ e_i=\nu_i,
\end{equation}
and, for $i\leq j$:
\begin{equation}\eqlabel{6.1.2}
e_j\circ e_i=e_i\circ e_j=e_i.
\end{equation}
\equref{6.1.1} is immediate; \equref{6.1.2} can be seen as follows:
\begin{eqnarray*}
e_j\circ e_i &=& e_j\circ \mu_i\circ \nu_i=e_j\circ \mu_j\circ\mu_{ji}\circ \nu_i\\
&\equal{\equref{6.1.1}}& \mu_j\circ\mu_{ji}\circ \nu_i=\mu_i\circ \nu_i=e_i;\\
e_i\circ e_j&=& \mu_i\circ \nu_i\circ e_j=\mu_i\circ \nu_{ij}\circ\nu_{j}\circ e_j\\
&\equal{\equref{6.1.1}}&\mu_i\circ \nu_{ij}\circ\nu_{j}=\mu_i\circ \nu_i=e_i.
\end{eqnarray*}

\begin{lemma}\lelabel{6.1}
Let $i\in I$ and $t\in T$. There exists $t_i\in T_i$ such that
$$t=\mu_i\circ t_i\circ \nu_i$$
if and only if
$$t=e_i\circ t\circ e_i.$$
In this situation,
$t_i$ is unique, and is given by the formula $t_i=\nu_i\circ t\circ \mu_i$; furthermore,
for every $j\geq i$, $t=\mu_j\circ t_j\circ \nu_j$, with
\begin{equation}\eqlabel{6.1.3}
t_j=\mu_{ji}\circ t_i\circ \nu_{ij}.
\end{equation}
\end{lemma}

\begin{proof}
We leave the first part as an easy exercise to the reader. For $j\geq i$, we compute
$$\mu_j\circ \mu_{ji}\circ t_i\circ \nu_{ij}\circ \nu_j=\mu_i\circ t_i\circ \nu_i=t.$$
\end{proof}

\begin{proposition}\prlabel{6.2}
$T^\dagger=\{t\in T~|~\exists i\in I:~t=e_i\circ t\circ e_i\}$ is a subalgebra of $T$
with idempotent local units. In particular, $T^\dagger$ is a firm $k$-algebra.
$\colim \ul{T}=(T^\dagger,\rho)$, with
$$\rho_i:\ T_i\to T^\dagger,~~\rho_i(t_i)=\mu_i\circ t_i\circ \nu_i.$$
\end{proposition}

\begin{proof}
It is clear that the $e_i$ form a set of
idempotent local units. It follows from \leref{6.1} that $T_i=e_iTe_i$.
$(T^\dagger,\rho)$ is a cocone on $T$ since, for all $j\geq i$
and $t_i\in T_i$, we have
$$(\rho_j\circ \rho_{ji})(t_i)=\mu_j\circ \mu_{ji}\circ t_i\circ\nu_{ij}\circ \nu_j
=\mu_i\circ  t_i\circ\nu_i=\rho_i(t_i).$$
Assume that $(M,m)$ is another cocone on $T^\dagger$. This means that
$m_i:\ T_i\to M$ and $m_j\circ \rho_{ji}=m_i$ if $i\leq j$. The map
$f:\ T^\dagger\to M$,
$f(\mu_i\circ t_i\circ \nu_i)=m_i(t_i)$,
is well-defined: assume that $t=\mu_i\circ t_i\circ \nu_i=\mu_j\circ t_j\circ \nu_j$.
Take $k\geq i,j$. By \equref{6.1.1}, $t=\mu_k\circ t_k\circ \nu_k$, with
$t_k=\mu_{ki}\circ t_i\circ \nu_{ik}=\rho_{ki}(t_i)$, and it follows that
$m_k(t_k)=m_k(\rho_{ki}(t_i))=m_i(t_i)$. In a similar way, we have that
$m_k(t_k)=m_j(t_j)$.\\
Finally, for every $i\in I$ and $t_i\in T_i$, we have that
$(f\circ \rho_i)(t_i)=f(\mu_i\circ t_i\circ\nu_i)=m_i(t_i)$.
\end{proof}

It is easy to show that $P_i\cong e_i\cdot P$, so that the comatrix coring
$P^{\dagger}\ot_{T^{\dagger}}P$ is a special case of the comatrix coring studied
in \seref{4}.
In general, $P^{\dagger}$ is a proper submodule of $P^*$ and
$T^\dagger$ is a proper subalgebra of $T$. But we have the following remarkable result.

\begin{proposition}\prlabel{6.3}
The map
$$\kappa:\ P^{\dagger}\ot_{T^{\dagger}}P\to
P^*\ot_TP,~~\kappa(\varphi\ot_{T^{\dagger}}p)=\varphi\ot_T p$$
is an isomorphism of $A$-bimodules.
\end{proposition}

\begin{proof}
We first define map
$\lambda:\ P^*\ot P\to P^{\dagger}\ot_{T^{\dagger}} P$
as follows: take $\varphi\in P^*$ and $p\in P$. There exists $i\in I$
such that $p=\sigma_i(p_i)$, and we define
$$\lambda(\varphi \ot p)=\varphi\circ\sigma_i\circ \tau_i\ot_{T^{\dagger}} p.$$
The right hand side does not depend on the choice of $i$: assume that $j\in I$ is such that
$p=\sigma_j(p_j)$ for some $p_j\in P_j$, and take $k\geq i,j$. Then we have that
$p=\sigma_k(p_k)$ with $p_k=\sigma_{ki}(p_i)$. We compute that
\begin{equation}\eqlabel{6.3.1}
\sigma_k\circ \tau_k\circ \sigma_i\circ \tau_i=
\sigma_k\circ \tau_k\circ \sigma_k\circ \sigma_{ki}\circ \tau_i
=\sigma_k\circ \sigma_{ki}\circ \tau_i= \sigma_i\circ \tau_i,
\end{equation}
hence
\begin{eqnarray*}
&&\hspace*{-2cm}
(\varphi\circ\sigma_i\circ\tau_i)\ot_{T^{\dagger}} p=
(\varphi\circ\sigma_k\circ \tau_k\circ\sigma_i\circ\tau_i)\ot_{T^{\dagger}} p\\
&=& (\varphi\circ\sigma_k\circ \tau_k)\ot_{T^{\dagger}} (\sigma_i\circ\tau_i)(p)
= (\varphi\circ\sigma_k\circ \tau_k)\ot_{T^{\dagger}} p,
\end{eqnarray*}
and, in a similar way,
$$(\varphi\circ\sigma_j\circ\tau_j)\ot_{T^{\dagger}} p=(\varphi\circ\sigma_k\circ \tau_k)\ot_{T^{\dagger}} p.$$
Our next aim is to show that
\begin{equation}\eqlabel{6.3.2}
\lambda(\varphi\ot t\cdot p)=\lambda(\varphi \cdot t\ot p),
\end{equation}
for all $\varphi\in P^*$, $t\in T$, and $p\in P$. There exists $i\in I$ such that
$$p=\sigma_i(\tau_i(p))=\sigma_i(p_i)~~{\rm and}~~
t\cdot p=(\sigma_i\circ\tau_i)(t\cdot p)=(\sigma_i\circ\tau_i\circ \omega(t)\circ\sigma_i)(p_i).$$
For all $k\geq i$, we then also have that
\begin{equation}\eqlabel{6.3.3}
p=\sigma_k(p_k)~~{\rm and}~~
t\cdot p=(\sigma_k\circ\tau_k\circ \omega(t)\circ\sigma_k)(p_k).
\end{equation}
For all $p\in P$, we have
$$\tau_i(p)=\sum z_iz_i^*(\tau_i(p)),$$
hence
$$(\omega(t)\circ\sigma_i\circ\tau_i)(p)=
\sum (\omega(t)\circ\sigma_i)(z_i)z_i^*(\tau_i(p)).$$
There exists $k\in I$ such that all $(\omega(t)\circ\sigma_i)(z_i)\in \sigma_k(P_k)$,
or $(\omega(t)\circ\sigma_i)(z_i)=(\sigma_k\circ\tau_k\circ \omega(t)\circ\sigma_i)(z_i).$
Then we find for all $p\in P$ that
\begin{eqnarray*}
&&\hspace*{-15mm}
(\omega(t)\circ\sigma_i\circ\tau_i)(p)=
\sum (\omega(t)\circ\sigma_i)(z_i)z_i^*(\tau_i(p))\\
&=&\sum (\sigma_k\circ\tau_k\circ\omega(t)\circ\sigma_i)(z_i)z_i^*(\tau_i(p))=
(\sigma_k\circ\tau_k\circ\omega(t)\circ\sigma_i\circ\tau_i)(p).
\end{eqnarray*}
We can take $k\geq i$. Using \equref{6.3.1}, we then find
\begin{equation}\eqlabel{6.3.4}
\omega(t)\circ\sigma_i\circ\tau_i=
\sigma_k\circ\tau_k\circ\omega(t)\circ\sigma_i\circ\tau_i=
\sigma_k\circ\tau_k\circ\omega(t)\circ\sigma_k\circ\tau_k\circ\sigma_i\circ\tau_i.
\end{equation}
We now compute
\begin{eqnarray*}
&&\hspace*{-2cm}
\lambda(\varphi t\ot p)=
\varphi\circ \omega(t)\circ\sigma_i\circ\tau_i\ot_{T^{\dagger}}
\sigma_i(p_i)\\
&\equal{\equref{6.3.4}}&
\varphi\circ \sigma_k\circ\tau_k\circ\omega(t)\circ\sigma_k\circ\tau_k\circ\sigma_i\circ\tau_i
\ot_{T^{\dagger}}\sigma_i(p_i)\\
&=&
\varphi\circ \sigma_k\circ\tau_k\circ\omega(t)\circ\sigma_k\circ\tau_k
\ot_{T^{\dagger}}(\sigma_i\circ\tau_i\circ\sigma_i)(p_i)\\
&=&
\varphi\circ \sigma_k\circ\tau_k\circ\omega(t)\circ\sigma_k\circ\tau_k
\ot_{T^{\dagger}}\sigma_i(p_i)\\
&=&
\varphi\circ \sigma_k\circ\tau_k\circ \sigma_k\circ\tau_k\circ\omega(t)\circ\sigma_k\circ\tau_k
\ot_{T^{\dagger}}\sigma_k(p_k)\\
&=&
\varphi\circ \sigma_k\circ\tau_k \ot_{T^{\dagger}}
(\sigma_k\circ\tau_k\circ\omega(t)\circ\sigma_k\circ\tau_k\circ\sigma_k)(p_k)\\
&\equal{\equref{6.3.3}}& \varphi\circ \sigma_k\circ\tau_k\ot_{T^{\dagger}}tp=\lambda(\varphi\ot tp),
\end{eqnarray*}
proving \equref{6.3.2}. We conclude that $\lambda$ induces a well-defined map
$\lambda:\ P^*\ot_TP\to P^{\dagger}\ot_{T^{\dagger}}P$.
Let us finally show that $\lambda$ is the inverse of $\kappa$. Take $\varphi\in P^{\dagger}$ and $p\in P$. Then there exists $i\in I$ such that $\varphi=\varphi_i\circ \tau_i$ and
$p=\sigma_i(p_i)$ for some $p_i\in P_i$, $\varphi_i\in P_i^*$. Then
$$
\lambda(\kappa(\varphi\ot_{T^{\dagger}}p))=\lambda(\varphi\ot_T p)
= \varphi\circ\sigma_i\circ\tau_i \ot_{T^{\dagger}}p
= \varphi\ot_{T^{\dagger}}(\sigma_i\circ\tau_i )(p)=\varphi\ot_{T^{\dagger}}p.
$$
Take $\varphi\in P^*$, and $p=\sigma_i(p_i)\in P$. Then
$$
\kappa(\lambda(\varphi\ot_T p))=\kappa(\varphi\circ\sigma_i\circ\tau_i
\ot_{T^{\dagger}}p)=\varphi\circ\sigma_i\circ\tau_i
\ot_{{T}}p=\varphi\ot_T (\sigma_i\circ\tau_i)(p)=\varphi\ot_Tp.$$
\end{proof}

We will now describe the infinite comatrix coring $P^{\dagger}\ot_{T^{\dagger}}P$
as the colimit of a richer system. On $I\times I$, we define a preorder as follows.
\begin{itemize}
\item $(i,i)\leq (j,j)$ if $i\leq j$ in $I$;
\item $(i,j)\leq (i,i)$, for all $i,j\in I$;
\item $(i,j)\leq (j,j)$, for all $i,j\in I$.
\end{itemize}
This preorder induces a partial order $\leq$ on $I\times I$. We have a corresponding
category $\Yy$. If $i\leq j$ in $I$, then the corresponding morphism $(i,i)\to (j,j)$
in $\Yy$ is denoted by $a_{ji}$. The morphism $(i,j)\to(i,i)$ is denoted by $l_{ij}$,
and the morphism $(i,j)\to(j,j)$ by $r_{ij}$. Note that we have a functor
$\xi:\ \Zz\to \Yy$, $\xi(i)=(i,i)$, $\xi(a_{ji})=a_{ji}$.

\begin{proposition}\prlabel{6.4}
We have a functor $\ul{F}:\ \Yy\to {}_A\Mm_A$ such that $\ul{F}\circ \xi=\ul{G}$.
\end{proposition}

\begin{proof}
For $i,j\in I$, $T_{ji}=\Hom_\Aa(M_i,M_j)$ is a $(T_j,T_i)$-bimodule, and we have
$$F(i,j)=P_j^*\ot_{T_j}T_{ji}\ot_{T_i}P_i.$$
We now define $\ul{F}$ on the morphisms. Let $\ul{F}(a_{ji})=\ul{G}(a_{ji})$; $\ul{F}(l_{ij})$
and $\ul{F}(r_{ij})$ are given by
$$\ul{F}(l_{ij}):\ P_j^*\ot_{T_j}T_{ji}\ot_{T_i}P_i\to P_i^*\ot_{T_i} P_i,~~
\ul{F}(l_{ij})(\varphi_j\ot_{T_j}t_{ji}\ot_{T_i}p_i)=\varphi_j\circ t_{ji}\ot_{T_i}p_i;$$
$$\ul{F}(r_{ij}):\ P_j^*\ot_{T_j}T_{ji}\ot_{T_i}P_i\to P_j^*\ot_{T_j} P_j,~~
\ul{F}(r_{ij})(\varphi_j\ot_{T_j}t_{ji}\ot_{T_i}p_i)=\varphi_j\ot_{T_j} t_{ji}(p_i).$$
We have to prove that $\ul{F}(a_{ji})\circ \ul{F}(l_{ij})=\ul{F}(r_{ij})$ if $i\leq j$. We compute
easily that
\begin{eqnarray*}
&&\hspace*{-2cm}
\ul{F}(a_{ji})\circ \ul{F}(l_{ij})(\varphi_j\ot_{T_j}t_{ji}\ot_{T_i}p_i)=
\varphi_j\circ t_{ji}\circ\beta(a_{ji})\ot_{T_j}\alpha(a_{ji})(p_i)\\
&=&\varphi_j\ot_{T_j}(t_{ji}\circ\beta(a_{ji})\circ \alpha(a_{ji}))(p_i)\\
&=& \varphi_j\ot_{T_j}t_{ji}(p_i)=\ul{F}(r_{ij})(\varphi_j\ot_{T_j}t_{ji}\ot_{T_i}p_i).
\end{eqnarray*}
In a similar way, we prove that
$\ul{F}(a_{ji})\circ \ul{F}(r_{ij})=\ul{F}(l_{ij})$ if $i\leq j$. All other verifications are easy.
\end{proof}

\begin{proposition}\prlabel{6.5}
$\colim \ul{F}= (P^{\dagger}\ot_{T^{\dagger}}P,f)$ with
$$f_{ij}=g_i\circ \ul{F}(l_{ij})=g_j\circ \ul{G}(r_{ij}),$$
for all $i,j\in I$.
\end{proposition}

\begin{proof}
It is easy to show that $(P^{\dagger}\ot_{T^{\dagger}}P,f)$ is a
cocone on $F$. If $(M,m)$ is another cocone on $F$, then we have a cocone $(M,n)$
on $G$, with $n_i=m_{(i,i)}$. We then have an $A$-bimodule map 
$f:\ P^{\dagger}\ot_{T^{\dagger}}P$, and it is straightforward to show that
it satisfies the necessary requirements.
\end{proof}

\section{Split direct systems of Galois comodules}\selabel{7}
Let $A$ be a $k$-algebra (with unit), and $\Cc$ an $A$-coring. By \cite[18.12]{BrzezinskiWisbauer} the category $\Mm^\Cc$ contains direct sums and cokernels. Consequently $\Mm^\Cc$ contains colimits, so in particular directed limits. Moreover, the forgetful
functor $\omega:\ \Mm^\Cc\to \Mm_A$ has a right adjoint, so it preserves
colimits (see for example \cite[Sec. V.5]{McLaine}). Hence we can apply the
results of \seref{6} in the situation where $\Aa=\Mm^\Cc$.

Now we consider a split direct system $\ul{M}^s:\ \Zz\to \Mm^\Cc_{\rm fgp}$,
$\ul{M}^s(i)=M_i$ and $\ul{M}^s(a_{ji})=(\mu_{ji},\nu_{ij})$. Here $\Mm^\Cc_{\rm fgp}$ denotes the category of right $\Cc$-comodules that are finitely generated and projective as right $A$-module.
We can compute $\colim \ul{M}=(M,\mu)$ in $\Mm^\Cc$, and from \prref{3.1} we know that there
exist left inverses $\nu_i$ of the $\mu_i$. As in \seref{6}, let
$T_i=\End^\Cc(M_i)$, and $(T^\dagger,t)=\colim T$. We have the associated
comatrix coring $\Gg= M^\dagger\ot_{T^\dagger} M$. For every $i\in I$, we have
a morphism of corings
$$\can_i:\ G(i)=M_i^*\ot_{T_i} M_i\to \Cc,~~
\can_i(\varphi_i\ot_{T_i}m_i)=\varphi_i(m_{i[0]})m_{i[1]}.$$

\begin{lemma}\lelabel{7.2}
$(\Cc,\can)$ is a cocone on $\ul{G}:\ \Zz\to {}_A\Mm_A$
\end{lemma}

\begin{proof}
For $i\leq j$ in $I$, $m_i\in M_i$ and $\varphi_i\in M_i^*$, we calculate that
\begin{eqnarray*}
&&\hspace*{-2cm}
(\can_j\circ\ul{G}(a_{ji}))(\varphi_i\ot_{T_i}m_i)=
\can_j(\varphi_i\circ\nu_{ij}\ot_{T_j}\mu_{ji}(m_i))\\
&=& (\varphi_i\circ\nu_{ij})(\mu_{ji}(m_i)_{[0]})\mu_{ji}(m_i)_{[1]}
=\varphi_i(\tau_{ij}(\mu_{ji}(m_{i[0]})))m_{i[1]}\\
&=&\varphi_i(m_{i[0]})m_{i[1]}=\can_i(\varphi_i\ot_{T_i}m_i),
\end{eqnarray*}
where we used the fact that $\mu_{ji}$ is right $\Cc$-colinear.
\end{proof}

\begin{proposition}\prlabel{7.3}
There exists a unique morphism of corings
$$\can:\ \Gg=M^\dagger\ot_{T^\dagger} M\to \Cc$$
such that
\begin{equation}\eqlabel{7.3.1}
\can(\varphi_i\circ\nu_i\ot_{T^\dagger}\mu_i(m_i))=\varphi_i(m_{i[0]})m_{i[1]},
\end{equation}
for all $i\in I$, $m_i\in M_i$ and $\varphi_i\in M_i^*$.
\end{proposition}

\begin{proof}
We have seen in \prref{4.9} that $\colim\ul{G}=(\Gg,g)$. It follows from \leref{7.2}
and the universal property of colimits that there exists a unique $(A,A)$-bimodule
map $\can:\ \Gg\to \Cc$ satisfying \equref{7.3.1}. The proof is finished if we can
show that $\can$ is a map of corings. As before, let $E_i=\sum z_i\ot_Az_i^*$
be a finite dual basis of $M$ as a right $A$-module. For all $m_i\in M_i$, we have that
$\sum z_iz_i^*(m_i)=m_i$.
Since $\rho_i$ is right $A$-linear, we have
$$\sum z_{i[0]}\ot_A z_{i[1]}z_i^*(m_i)=m_{i[0]}\ot_A m_{i[1]},$$
hence, for all $\varphi_i\in M_i$:
\begin{equation}\eqlabel{7.3.2}
\sum \varphi_i(z_{i[0]})z_{i[1]}z_i^*(m_i)=\varphi_i(m_{i[0]}) m_{i[1]}.
\end{equation}
Then we compute
\begin{eqnarray*}
&&\hspace*{-2cm}
(\can\ot_A\can)\Delta(\varphi_i\circ\nu_i\ot_{T^\dagger}\mu_i(m_i))
= \sum \varphi_i(z_{i[0]})z_{i[1]}\ot_A z_i^*(p_{i[0]})p_{i[1]}\\
&=& \sum \varphi_i(z_{i[0]})z_{i[1]}z_i^*(p_{i[0]})\ot_A p_{i[1]}
\equal{\equref{7.3.2}} \varphi_i(p_{i[0]})p_{i[1]}\ot_A p_{i[2]}\\
&=&\Delta(\varphi_i(p_{i[0]})p_{i[1]})=\Delta(\can(\varphi_i\circ\nu_i\ot_{T^\dagger}\mu_i(m_i))).
\end{eqnarray*}
Finally
\begin{eqnarray*}
&&\hspace*{-2cm}
\varepsilon(\can(\varphi_i\circ\nu_i\ot_{T^\dagger}\mu_i(m_i)))
=\varepsilon(\varphi_i(m_{i[0]})m_{i[1]})
= \varphi_i(m_{i[0]})\varepsilon(m_{i[1]})\\
&=&\varphi_i(m_{i[0]}\varepsilon(m_{i[1]}))
= \varphi_i(m_i)=\varepsilon(\varphi_i\circ\nu_i\ot_{T^\dagger}\mu_i(m_i)).
\end{eqnarray*}
\end{proof}

We call $\ul{M}^s$ a {\sl split direct system of Galois $\Cc$-comodules}
if $\can:\ \Gg\to \Dd$ is an isomorphism of corings. In this situation, we have that
the categories $\Mm^\Gg$ and $\Mm^\Cc$ are isomorphic. From \thref{1.3}, we then
immediately obtain the following result:

\begin{theorem}\thlabel{7.4}
Let $\Cc$ be an $A$-coring, and 
$\ul{M}^s$ a split direct system of Galois $\Cc$-comodules. Let $\colim \ul{M}=
(M,\mu)$, $\colim \End^\Cc(\ul{M})=(T^\dagger, \rho)$ and $\Gg=M^\dagger\ot_{T^\dagger}
M$ the associated comatrix coring. If $M$ is faithfully flat as a left $T^\dagger$-module,
then the categories $\Mm_{T\dagger}$ and $\Mm^\Cc$ are equivalent.
\end{theorem}

\begin{center}
{\sc Acknowledgements}
\end{center}
We thank the referee for his detailed comments.

\end{document}